\documentclass[11pt, a4paper]{article}
\usepackage{amsmath}
\usepackage{wasysym}
\usepackage{latexsym}
\usepackage{amsfonts}
\usepackage{mathrsfs}
\usepackage{amssymb}
\usepackage{ifsym}
\usepackage{dsfont}
\usepackage[all]{xy}
\usepackage{authblk}
\usepackage{float}

\usepackage{siunitx} 

\newtheorem{Theorems1}{Theorem}[section]

\newtheorem{Coroll1}[Theorems1]{Corollary}
\newtheorem{Lemma1}[Theorems1]{Lemma}

\newtheorem{Examp1}[Theorems1]{Example}

\newtheorem{Quest1}[Theorems1]{Question}

\newtheorem{Definitions1}[Theorems1]{Definition}

\newenvironment{proof}[1][Proof]{\begin{trivlist}
\item[\hskip \labelsep {\bfseries #1}]}{\end{trivlist}}

\newcommand{\qed}{\nobreak \ifvmode \relax \else
      \ifdim\lastskip<1.5em \hskip-\lastskip
      \hskip1.5em plus0em minus0.5em \fi \nobreak
      \vrule height0.75em width0.5em depth0.25em\fi}

\begin{document}
\title{End extending models of set theory via \break power admissible covers}

\author[1]{Zachiri McKenzie}
\author[2]{Ali Enayat}
\affil[1]{{\tt zach.mckenzie@gmail.com}}
\affil[2]{{\tt ali.enayat@gu.se}}
\maketitle

\begin{abstract}
Motivated by problems involving end extensions of models of set theory, we develop the rudiments of the \emph{power admissible cover} construction (over ill-founded models of set theory), an extension of the machinery of admissible covers invented by Barwise as a versatile tool for generalising model-theoretic results about countable well-founded models of set theory to countable ill-founded ones. Our development of the power admissible machinery allows us to obtain new results concerning \emph{powerset-preserving end extensions and rank extensions} of countable models of subsystems of $\mathsf{ZFC}$. The canonical extension $\mathsf{KP}^\mathcal{P}$ of Kripke-Platek set theory $\mathsf{KP}$ plays a key role in our work; one of our results refines a theorem of Rathjen by showing that $\Sigma_1^\mathcal{P}$-\textsf{Foundation} is provable in $\mathsf{KP}^\mathcal{P}$ (without invoking the axiom of choice).
\end{abstract}

\section[Introduction]{Introduction}

The \emph{admissible cover} machinery was introduced by Barwise in the Appendix of his venerable book \cite{bar75} on admissible set theory. Admissible covers allow one to extend the range of infinitary compactness arguments from the domain of countable well-founded models of $\mathsf{KP}$ (Kripke-Platek set theory) to countable \emph{ill-founded} models of  $\mathsf{KP}$. For example, Barwise uses admissible covers in his book to prove a striking result:  \emph{Every countable model of $\mathsf{ZF}$ has an extension to a model of $\mathsf{ZF+V=L}$}.\footnote{This end extension result, together with certain elaborations of it, first appeared in an earlier paper of Barwise \cite{bar71}. It is also noteworthy that, as shown recently by Hamkins \cite{ham18}, Barwise's end extension theorem can also be proved using more classical techniques (without appealing to methods of admissible set theory).}  Admissible covers also appear in the work of Ressayre \cite{res87}, who showed that the results presented in the Appendix of \cite{bar75} pertaining to $\mathsf{KP}$ do not depend on the availability of the full scheme of foundation among the axioms of $\mathsf{KP}$; more specifically, they only require the scheme of foundation for $\Sigma_1 \cup \Pi_1$-formulae.\footnote{Note that the full scheme of foundation is included in the axioms of $\mathsf{KP}$ in Barwise's treatment \cite{bar75}. However, we follow the convention proposed by Mathias to only include $\Pi_1$-\textsf{Foundation} in the axiomatisation of $\mathsf{KP}$; this is informed by the fact, demonstrated by Mathias \cite{mat01}, that many (but not all) results about Barwise's  $\mathsf{KP}$ can be carried out within Mathias' $\mathsf{KP}$.} Admissible covers were used more recently by Williams \cite{wil18}, to show that certain class theories (including Kelley-Morse class theory) fail to have minimum transitive models (this result of Williams also appears in their paper \cite{wil19}, but with a different proof). In this paper we explore the variant \emph{power admissible cover} of the notion of admissible cover in order to obtain new results in the model theory of set theory. The main inspiration for our results on end extensions arose from our joint work with Kaufmann \cite{ekm18} on automorphisms of models of set theory (see Theorem \ref{thm: aut}).

The canonical extension $\mathsf{KP}^\mathcal{P}$ of Kripke-Platek set theory $\mathsf{KP}$ plays a key role in our work. $\mathsf{KP}^\mathcal{P}$ is intimately related to Friedman's so-called power admissble system $\mathsf{PAdm^{s}}$, whose well-founded models are the the so-called \emph{power admissible sets}  \cite{fri73}.\footnote{The precise relationship between  Friedman's system and $\mathsf{KP}^\mathcal{P}$ is worked out in Section 6.19 of \cite{mat01}.}
These two systems can accommodate constructions by $\Sigma_1$-recursions relative to the power set operation. The system $\mathsf{KP}^\mathcal{P}$ has been closely studied by Mathias \cite{mat01} and Rathjen \cite{rat14}, \cite{rat20}. In the latter paper Rathjen proves that $\Sigma_1^\mathcal{P}$-\textsf{Foundation} is provable in $\mathsf{KP}^\mathcal{P} + \mathsf{AC}$ (where $\mathsf{AC}$ is the axiom of choice).

The highlights of the paper are as follows. In Corollary 3.3 we refine Rathjen's aforementioned result by showing that $\Sigma_1^\mathcal{P}$-\textsf{Foundation} is provable outright in $\mathsf{KP}^\mathcal{P}$. The rudiments of power admissible covers are developed in Section 4. In Section 5 the machinery of power admissible covers is put together with results of earlier sections to establish new results about powerset-preserving end extensions and rank extensions of models of set theory.  For example in Theorem 5.7 we show that every countable model of $\mathcal{M}\models\mathsf{KP}^\mathcal{P}$ has a topless rank extension, i.e., $\mathcal{M}$ has a proper rank extension $\mathcal{N}\models\mathsf{KP}^\mathcal{P}$ such that $\mathsf{Ord}^{\mathcal{N}} \setminus \mathsf{Ord}^{\mathcal{M} }$ has no least element. This result generalises a classical theorem of Friedman that shows that every countable \emph{well-founded} model of $\mathsf{KP}^\mathcal{P}$ has a topless rank extension.

\section[Background]{Background}

We use $\mathcal{L}$ throughout the paper to denote the \emph{language $\{\in, =\}$ of set theory}. We will make reference to generalisations of the L\'{e}vy hierarchy of formulae in languages extending $\mathcal{L}$ that possibly contain constant and function symbols.

Let $\mathcal{L}^\prime$ be a language extending $\mathcal{L}$. We use $\Delta_0(\mathcal{L}^\prime)$ to denote the smallest class of $\mathcal{L}^\prime$-formulae that is closed under the connectives of propositional logic and quantification in the form $\exists x \in t$ and $\forall x \in t$, where $t$ is a term of $\mathcal{L}^\prime$ and $x$ is a variable that does not appear in $t$. The classes $\Sigma_1(\mathcal{L}^\prime)$, $\Pi_1(\mathcal{L}^\prime)$, $\Sigma_2(\mathcal{L}^\prime)$, \ldots are defined inductively from $\Delta_0(\mathcal{L}^\prime)$ in the usual way. We will write $\Delta_0$, $\Sigma_1$, $\Pi_1$, \ldots instead of $\Delta_0(\mathcal{L})$, $\Sigma_1(\mathcal{L})$, $\Pi_1(\mathcal{L})$, \ldots, and we will use $\Pi_\infty$ and $\Pi_\infty(\mathcal{L}^\prime)$ to denote the class of all $\mathcal{L}$-formulae and $\mathcal{L}^\prime$-formulae respectively. An $\mathcal{L}^\prime$-formula is $\Delta_n(\mathcal{L}^\prime)$, for $n > 0$, if it is equivalent to both a $\Sigma_n(\mathcal{L}^\prime)$-formula and a $\Pi_n(\mathcal{L}^\prime)$-formula.

The class $\Delta_0^\mathcal{P}$ is the smallest class of $\mathcal{L}$-formulae that is closed under the connectives of propositional logic and quantification in the form $\mathcal{Q} x \subseteq y$ and $\mathcal{Q} x \in y$ where $\mathcal{Q}$ is $\exists$ or $\forall$, and $x$ and $y$ are distinct variables. The \emph{Takahashi classes} $\Delta_1^\mathcal{P}$, $\Sigma_1^\mathcal{P}$, $\Pi_n^\mathcal{P}$, \ldots are defined from $\Delta_0^\mathcal{P}$ in the same way as the classes $\Delta_1$, $\Sigma_1$, $\Pi_1$, \ldots are defined from $\Delta_0$. If $\Gamma$ is a collection of $\mathcal{L}^\prime$-formulae and $T$ is an $\mathcal{L}^\prime$-theory, then we write $\Gamma^T$ for the class of $\mathcal{L}^\prime$-formulae that are provably in $T$ equivalent to a formula in $\Gamma$.


We will use capital calligraphic font letters ($\mathcal{M}$, $\mathcal{N}$, \ldots) to denote $\mathcal{L}$-structures. If $\mathcal{M}$ is an $\mathcal{L}$-structure, then, unless we explicitly state otherwise, $M$ will be used to denote the underlying set of $\mathcal{M}$ and $\mathsf{E}^\mathcal{M}$ will be used to denote the interpretation of $\in$ in $\mathcal{M}$.

Let $\mathcal{L}^\prime$ be a language extending $\mathcal{L}$ and let $\mathcal{M}$ be an $\mathcal{L}^\prime$-structure with underlying set $M$. If $a \in M$, then $a^*$ is defined as follows:
$$a^*:=\{x \in M \mid \mathcal{M} \models (x \in a)\},$$ 
as long as the structure $\mathcal{M}$ is clear from the context. Let $\Gamma$ be a class of formulae. We say that $A \subseteq M$ is {\it $\Gamma$-definable over $\mathcal{M}$} if there exists a $\Gamma$-formula $\phi(x, \vec{z})$ and $\vec{a} \in M$ such that $A=\{x \in M \mid \mathcal{M} \models \phi(x, \vec{a})\}$.

Let $\mathcal{M}$ and $\mathcal{N}$ be $\mathcal{L}^\prime$-structures. We will partake in the common abuse of notation and write $\mathcal{M} \subseteq \mathcal{N}$ if $\mathcal{M}$ is a substructure of $\mathcal{N}$.
\begin{itemize}
\item We say that $\mathcal{N}$ is an \emph{end extension} of $\mathcal{M}$, and write $\mathcal{M} \subseteq_e \mathcal{N}$, if $\mathcal{M} \subseteq \mathcal{N}$ and for all $x, y \in N$, if $y \in M$ and $\mathcal{N} \models (x \in y)$, then $x \in M$.
\item We say that $\mathcal{N}$ is a \emph{powerset-preserving end extension} of $\mathcal{M}$, and write $\mathcal{M} \subseteq_e^\mathcal{P} \mathcal{N}$, if $\mathcal{M} \subseteq_e \mathcal{N}$ and for all $x, y \in N$, if $y \in M$ and $\mathcal{N} \models (x \subseteq y)$, then $x \in M$.
\item We say that $\mathcal{N}$ is a \emph{topless powerset-preserving end extension} of $\mathcal{M}$, and write $\mathcal{M} \subseteq_{\mathsf{topless}}^\mathcal{P} \mathcal{N}$, if $\mathcal{M} \subseteq_e^\mathcal{P} \mathcal{N}$, $M \neq N$ and for all $c \in N$, if $c^* \subseteq M$, then $c \in M$.
\item We say that $\mathcal{N}$ is a \emph{blunt powerset-preserving end extension} of $\mathcal{M}$, and write $\mathcal{M} \subseteq_{\mathsf{blunt}}^\mathcal{P} \mathcal{N}$, if $\mathcal{M} \subseteq_e^\mathcal{P} \mathcal{N}$, $M \neq N$ and $\mathcal{N}$ is not a topless powerset-preserving end extension of $\mathcal{M}$.
\end{itemize}

Let $\Gamma$ be a class of $\mathcal{L}$-formulae. The following define the restriction of the \textsf{ZF}-provable schemes  \textsf{Separation}, \textsf{Collection}, and \textsf{Foundation} to formulae in the class $\Gamma$:
\begin{itemize}
\item[]($\Gamma$-\textsf{Separation}) For all $\phi(x, \vec{z}) \in \Gamma$,
$$\forall \vec{z} \forall w \exists y \forall x(x \in y \iff (x \in w) \land \phi(x, \vec{z})).$$
\item[]($\Gamma$-\textsf{Collection}) For all $\phi(x, y, \vec{z}) \in \Gamma$,
$$\forall \vec{z} \forall w((\forall x \in w) \exists y \phi(x, y, \vec{z}) \Rightarrow \exists C (\forall x \in w)(\exists y \in C) \phi(x, y, \vec{z})).$$
\item[]($\Gamma$-\textsf{Foundation}) For all $\phi(x, \vec{z}) \in \Gamma$,
$$\forall \vec{z}(\exists x \phi(x, \vec{z}) \Rightarrow \exists y(\phi(y, \vec{z}) \land (\forall x \in y) \neg \phi(x, \vec{z}))).$$
If $\Gamma= \{x \in z\}$ then we will refer to $\Gamma$-\textsf{Foundation} as \textsf{Set-foundation}.
\end{itemize}
We will also make reference to the following fragments of \textsf{Separation} and \textsf{Foundation} for formulae that are $\Delta_n$ with parameters:
\begin{itemize}
\item[]($\Delta_n$-\textsf{Separation}) For all $\Sigma_n$-formulae, $\phi(x, \vec{z})$, and for all $\Pi_n$-formulae, $\psi(x, \vec{z})$,
$$\forall \vec{z}( \forall x (\phi(x, \vec{z}) \iff \psi(x, \vec{z})) \Rightarrow \forall w \exists y \forall x(x \in y \iff (x \in w) \land \phi(x, \vec{z}))).$$
\item[]($\Delta_n$-\textsf{Foundation}) For all $\Sigma_n$-formulae, $\phi(x, \vec{z})$, and for all $\Pi_n$-formulae, $\psi(x, \vec{z})$,
$$\forall \vec{z}(\forall x(\phi(x, \vec{z}) \iff \psi(x, \vec{z})) \Rightarrow (\exists x \phi(x, \vec{z}) \Rightarrow \exists y (\phi(y, \vec{z}) \land (\forall x \in y) \neg \phi(x, \vec{z}))).$$
\end{itemize}
Similar definitions can also be used to express $\Delta_n^\mathcal{P}$-\textsf{Separation} and $\Delta_n^\mathcal{P}$-\textsf{Foundation}.

We use $\mathsf{TCo}$ to denote the axiom that asserts that every set is contained in a transitive set.

We will consider extensions of the following subsystems of $\mathsf{ZFC}$:
\begin{itemize}
\item $\mathsf{S}_1$ is the $\mathcal{L}$-theory with axioms: \textsf{Extensionality}, \textsf{Emptyset}, \textsf{Pair}, \textsf{Union}, \textsf{Set difference}, and \textsf{Powerset}.
\item $\mathsf{M}$ is obtained from $\mathsf{S}_1$ by adding $\mathsf{TCo}$, \textsf{Infinity}, $\Delta_0$-\textsf{Separation}, and \textsf{Set-foundation}.
\item $\mathsf{Mac}$ is obtained from $\mathsf{M}$ by adding $\mathsf{AC}$ (the axiom of choice).
\item $\mathsf{M}^-$ is obtained from $\mathsf{M}$ by removing \textsf{Powerset}.
\item $\mathsf{KP}$ is the $\mathcal{L}$-theory with axioms: \textsf{Extensionality}, \textsf{Pair}, \textsf{Union}, $\Delta_0$-\textsf{Separation}, $\Delta_0$-\textsf{Collectio}n and $\Pi_1$-\textsf{Foundation}.
\item $\mathsf{KP}^-$ is obtained from $\mathsf{KP}$ by removing $\Pi_1$-\textsf{Foundation}.
\item $\mathsf{KPI}$ is obtained $\mathsf{KP}$ by adding \textsf{Infinity}.
\item $\mathsf{KP}^\mathcal{P}$ is obtained from $\mathsf{M}$ by adding $\Delta_0^\mathcal{P}$-\textsf{Collection} and $\Pi_1^\mathcal{P}$-\textsf{Foundation}.
\item $\mathsf{MOST}$ is obtained from $\mathsf{M}$ by adding $\Sigma_1$-\textsf{Separation} and $\mathsf{AC}$.
\end{itemize}
In subsystems of $\mathsf{ZFC}$ that include \textsf{Infinity} we can also consider the following restriction of $\Gamma$-\textsf{Foundation}:
\begin{itemize}
\item[]($\Gamma$-\textsf{Foundation on} $\omega$) For all $\phi(x, \vec{z}) \in \Gamma$,
$$\forall \vec{z}((\exists x \in \omega)\phi(x, \vec{z}) \Rightarrow (\exists y \in \omega)(\phi(y, \vec{z}) \land (\forall x \in y) \neg \phi(x, \vec{z}))).$$
\end{itemize}

The second family of theories that we will be concerned with are extensions of the variant of Kripke-Platek Set Theory with \emph{urelements} that is introduced in \cite[Appendix]{bar75}.

Let $\mathcal{L}^*$ be obtained from $\mathcal{L}$ by adding a second binary relation $\mathsf{E}$, a unary predicate $\mathsf{U}$, and a unary function symbol $\mathsf{F}$. The intended interpretation of $\mathsf{U}$ is to distinguish urelements from sets. The binary relation $\mathsf{E}$ is intended to be a membership relation that holds between urelements, and $\in$ is intended to be a membership relation that can hold between sets or urelements and sets.

Let $\mathcal{L}_\mathsf{P}^*$ be obtained from $\mathcal{L}^*$ by adding a new unary function symbol $\mathsf{P}$. An $\mathcal{L}_{\mathsf{P}}^{*}$-structure is a structure $\mathfrak{A}_\mathcal{M}= \langle \mathcal{M}; A, \in^\mathfrak{A}, \mathsf{F}^\mathfrak{A}, \mathsf{P}^\mathfrak{A} \rangle$, where $\mathcal{M}= \langle M, \mathsf{E}^\mathfrak{A} \rangle$, $M$ is the extension of $\mathsf{U}$, $A$ is the extension of $\neg \mathsf{U}$, $\in^\mathfrak{A}$ is the interpretation of $\in$, $\mathsf{E}^\mathfrak{A}$ is the interpretation of $\mathsf{E}$, $\mathsf{F}^\mathfrak{A}$ is the interpretation of $\mathsf{F}$, and $\mathsf{P}^\mathfrak{A}$ is the interpretation of $\mathsf{P}$.

$\mathcal{L}^*$-structures will be presented in the same format, but without an interpretation of $\mathsf{P}$. The $\mathcal{L}^*$- and $\mathcal{L}_{\mathsf{P}}^{*}$-theories presented below will ensure that $\mathsf{E}^\mathfrak{A} \subseteq M \times M$, and $\in^\mathfrak{A}~ \subseteq (M \cup A) \times A$.

Following \cite{bar75}, we simplify the presentation of $\mathcal{L}^*$- and $\mathcal{L}_{\mathsf{P}}^{*}$-formulae by treating these languages as two-sorted rather than one-sorted.

\medskip

When writing $\mathcal{L}^*$- and $\mathcal{L}_{\mathsf{P}}^{*}$-formulae we will use the convention below of Barwise \cite{bar75}.

\begin{itemize}

\item The variables $p, q, p_1, \ldots$ range over elements of the domain that satisfy $\mathsf{U}$ (urelements).

\item the variables $a, b, c, d, f, \ldots$ range over elements of the domain that satisfy $\neg \mathsf{U}$ (sets); and

\item the variables $x, y, z, w \ldots$ range over all elements of the domain.

\end{itemize}

\noindent Therefore, $\forall p (\cdots)$ is an abbreviation of $\forall x (\mathsf{U}(x) \Rightarrow \cdots)$, $\exists a(\cdots)$ is an abbreviation of $\exists x(\neg \mathsf{U}(x) \land \cdots)$, etc.

In section \ref{Sec:Cover}, we will see that certain $\mathcal{L}$-structures can interpret $\mathcal{L}^*$- and $\mathcal{L}_P^*$-structures in which the urelements are isomorphic to the original $\mathcal{L}$-structure. It is this interaction that motivates our unorthodox convention of using $\mathsf{E}^\mathcal{M}$, $\mathsf{E}^\mathcal{N}$, \ldots to denote the interpretation of $\in$ in the $\mathcal{L}$-structures $\mathcal{M}$, $\mathcal{N}$, \ldots It should be noted that this convention differs from Barwise \cite{bar75} where $E$ is consistently used to denote the interpretation of $\in$ in $\mathcal{L}$-structures.

The following are analogues of axioms, fragments of axiom schemes and fragments of theorem schemes of $\mathsf{ZFC}$ in the languages $\mathcal{L}^*$ and $\mathcal{L}_{\mathsf{P}}^{*}$:
\begin{itemize}
\item[](\textsf{Extensionality for sets}) $\forall a \forall b (a=b \iff \forall x (x \in a \iff x \in b))$.
\item[](\textsf{Pair}) $\forall x \forall y \exists a \forall z(z \in a \iff z=x \lor z=y)$.
\item[](\textsf{Union}) $\forall a \exists b (\forall y \in a)(\forall x \in y)(x \in b)$.
\end{itemize}
Let $\Gamma$ be a class of $\mathcal{L}_{\mathsf{P}}^{*}$-formulae.
\begin{itemize}
\item[]($\Gamma$-\textsf{Separation}) For all $\phi(x, \vec{z}) \in \Gamma$,
$$\forall \vec{z} \forall a \exists b \forall x(x \in b \iff (x \in a) \land \phi(x, \vec{z})).$$
\item[]($\Gamma$-\textsf{Collection}) For all $\phi(x, y, \vec{z}) \in \Gamma$,
$$\forall \vec{z} \forall a ((\forall x \in a)\exists y \phi(x, y, \vec{z}) \Rightarrow \exists b (\forall x \in a)(\exists y \in b) \phi(x, y, \vec{z})).$$
\item[]($\Gamma$-\textsf{Foundation}) For all $\phi(x, \vec{z}) \in \Gamma$,
$$\forall \vec{z} (\exists x \phi(x, \vec{z}) \Rightarrow \exists y(\phi(y, \vec{z}) \land (\forall w \in y) \neg \phi(w, \vec{z}))).$$
\end{itemize}
The following axiom in the language $\mathcal{L}^*$ describes the desired behaviour of the function symbol $\mathsf{F}$:
\begin{itemize}
\item[]($\dagger$) $\forall p \forall x(x \mathsf{E} p \iff x \in \mathsf{F}(p)) \land \forall a(\mathsf{F}(a)= \emptyset)$.
\end{itemize}
The next axiom, in the language $\mathcal{L}_{\mathsf{P}}^{*}$, says that the function symbol $\mathsf{P}$ is the usual powerset function:
\begin{itemize}
\item[]($\mathsf{Powerset}$) $\forall a \forall b(b \in \mathsf{P}(a) \iff b \subseteq a)$.
\end{itemize}
We will have cause to consider the following theories:
\begin{itemize}
\item $\mathsf{KPU}_{\mathbb{C}\mathrm{ov}}$ is the $\mathcal{L}^*$-theory with axioms: $\exists a(a=a)$, $\forall p \forall x(x \notin p)$, \textsf{Extensionality for sets}, \textsf{Pair}, \textsf{Union}, $\Delta_0(\mathcal{L}^*)$-\textsf{Separation}, $\Delta_0(\mathcal{L}^*)$-\textsf{Collection}, $\Pi_1(\mathcal{L}^*)$-\textsf{Foundation} and ($\dagger$).
\item $\mathsf{KPU}^{\mathcal{P}}_{\mathbb{C}\mathrm{ov}}$ is the $\mathcal{L}_{\mathsf{P}}^{*}$-theory obtained from $\mathsf{KPU}_{\mathbb{C}\mathrm{ov}}$ by adding \textsf{Powerset}, $\Delta_0(\mathcal{L}_{\mathsf{P}}^{*})$-\textsf{Separation}, $\Delta_0(\mathcal{L}_{\mathsf{P}}^{*})$-\textsf{Collection} and $\Pi_1(\mathcal{L}_{\mathsf{P}}^{*})$-\textsf{Foundation}.
\end{itemize}

\begin{Definitions1}
Let $\mathcal{M}= \langle M, \mathsf{E}^\mathcal{M} \rangle$ be an $\mathcal{L}$-structure.

An \textbf {admissible set covering} $\mathcal{M}$ is an $\mathcal{L}^*$-structure $$\mathfrak{A}_\mathcal{M}= \langle \mathcal{M}; A, \in^\mathfrak{A}, \mathsf{F}^\mathfrak{A} \rangle \models \mathsf{KPU}_{\mathbb{C}\mathrm{ov}}.$$  such that $\in^\mathfrak{A}$ is well-founded.

A \textbf{power admissible set covering} $\mathcal{M}$ is an $\mathcal{L}_{\mathsf{P}}^{*}$-structure $$\mathfrak{A}_\mathcal{M}= \langle \mathcal{M}; A, \in^\mathfrak{A}, \mathsf{F}^\mathfrak{A}, \mathsf{P}^\mathfrak{A} \rangle\models \mathsf{KPU}^{\mathcal{P}}_{\mathbb{C}\mathrm{ov}}$$ such that $\in^\mathfrak{A}$ is well-founded.

We use $\mathbb{C}\mathrm{ov}_\mathcal{M}= \langle \mathcal{M}; A_\mathcal{M}, \in, \mathsf{F}_\mathcal{M} \rangle$ to denote the smallest admissible set covering $\mathcal{M}$ whose membership relation $\in$ coincides with the membership relation of the metatheory.

We use $\mathbb{C}\mathrm{ov}_\mathcal{M}^\mathsf{P}=\langle \mathcal{M}; A_\mathcal{M}, \in, \mathsf{F}_\mathcal{M}, \mathsf{P}_\mathcal{M} \rangle$ to denote the smallest power admissible set covering $\mathcal{M}$ whose membership relation coincides with the membership relation of the metatheory.

Note that if $\mathfrak{A}_\mathcal{M}= \langle \mathcal{M}; A, \in^\mathfrak{A}, \mathsf{F}^\mathfrak{A}, \ldots \rangle$ is an admissible set covering $\mathcal{M}$, then $\mathfrak{A}_\mathcal{M}$ is isomorphic to a structure whose membership relation $\in$ is the membership relation of the metatheory.
\end{Definitions1}

\begin{Definitions1}
Let $\mathcal{M}= \langle M, \mathsf{E}^\mathcal{M} \rangle$ be an $\mathcal{L}$-structure, and let $$\mathfrak{A}_\mathcal{M}= \langle \mathcal{M}; A, \in^\mathfrak{A}, \mathsf{F}^\mathfrak{A}, \mathsf{P}^\mathfrak{A} \rangle \models \mathsf{KPU}^{\mathcal{P}}_{\mathbb{C}\mathrm{ov}}.$$

We use $\mathrm{WF}(A)$ to denote the largest $B \subseteq_e A$ such that $\langle B, \in^\mathfrak{A} \rangle$ is well-founded.

The \textbf{well-founded part} of $\mathfrak{A}_\mathcal{M}$ is the $\mathcal{L}_{\mathsf{P}}^{*}$-structure $$\mathrm{WF}(\mathfrak{A}_\mathcal{M})= \langle \mathcal{M}; \mathrm{WF}(A), \in^\mathfrak{A}, \mathsf{F}^\mathfrak{A}, \mathsf{P}^\mathfrak{A} \rangle.$$

Note that $\mathrm{WF}(\mathfrak{A}_\mathcal{M})$ is always isomorphic to an $\mathcal{L}_{\mathsf{P}}^{*}$-structure whose membership relation $\in$ coincides with the membership relation of the metatheory.
\end{Definitions1}

As usual, in the theories $\mathsf{M}^-$, $\mathsf{KP}^-$ and $\mathsf{KPU}_{\mathbb{C}\mathrm{ov}}$ the ordered pair $\langle x, y \rangle$ is coded by the set $\{\{x\}, \{x, y\}\}$. This definition ensures that there is a $\Delta_0$-formula $\mathsf{OP}(x)$ that says that $x$ is an ordered pair, and functions
$$\mathsf{fst}(\langle x, y \rangle)= x \textrm{ and } \mathsf{snd}(\langle x, y\rangle)= y,$$
whose graphs are defined by $\Delta_0$-formulae. In $\mathsf{KPU}_{\mathbb{C}\mathrm{ov}}$ the \emph{rank function}, $\rho$, and \emph{support function}, $\mathsf{sp}$, are defined by recursion:
$$\rho(p)= 0 \textrm{ for all urelements } p, \textrm{ and } \rho(a)= \sup\{ \rho(x)+1 \mid x \in a\} \textrm{ for all sets } a;$$
$$\mathsf{sp}(p)= \{p\} \textrm{ for all urelements } p, \textrm{ and } \mathsf{sp}(a)= \bigcup_{x \in a} \mathsf{sp}(x) \textrm{ for all sets } a.$$
The theory $\mathsf{KPU}_{\mathbb{C}\mathrm{ov}}$ proves that both of these are total and their graphs are $\Delta_1(\mathcal{L}^*)$. In the theory $\mathsf{KP}$, in which everything is a set, the rank function, $\rho$, is $\Delta_1$ and remains provably total. We say that $x$ is a \emph{pure set} if $\mathsf{sp}(x)= \emptyset$. We say that $x$ is an \emph{ordinal} if $x$ is a hereditarily transitive pure set; where:
$$\mathsf{Transitive}(x) \iff \neg \mathsf{U}(x) \land (\forall y \in x)(\forall z \in y)(z \in x), \textrm{ and}$$
$$\mathsf{Ord}(x) \iff (\mathsf{Transitive}(x) \land (\forall y \in x)(\mathsf{Transitive}(y)).$$
Therefore, both `$x$ is transitive' and `$x$ is an ordinal' can be expressed using $\Delta_0(\mathcal{L}^*)$-formulae. In the theories $\mathsf{M}^-$ and $\mathsf{KP}$, we can omit the reference to the predicate $\mathsf{U}$ in the definition of `$x$ is transitive', thus making both the property of being transitive and the property of being an ordinal into $\Delta_0$ properties.

The rank function allows us to strengthen the notion of powerset-preserving end extensions for models of $\mathsf{KP}$. Let $\mathcal{L}^\prime$ be a language extending $\mathcal{L}$. Let $\mathcal{M}$ and $\mathcal{N}$ be $\mathcal{L}^\prime$-structures that satisfy $\mathsf{KP}$.
\begin{itemize}
\item We say that $\mathcal{N}$ is a \emph{rank extension} of $\mathcal{M}$, and write $\mathcal{M} \subseteq_e^{\mathsf{rk}} \mathcal{N}$, if $\mathcal{M} \subseteq_e^\mathcal{P} \mathcal{N}$ and for all $x, y \in N$, if $y \in M$ and $\mathcal{N} \models (\rho(x) \leq \rho(y))$, then $x \in M$.
\item We say that $\mathcal{N}$ is a \emph{topless rank extension} of $\mathcal{M}$, and write $\mathcal{M} \subseteq_{\mathsf{topless}}^{\mathrm{rk}} \mathcal{N}$, if $\mathcal{M} \subseteq_e^{\mathrm{rk}} \mathcal{N}$, $M \neq N$ and for all $c \in N$, if $c^* \subseteq M$, then $c \in M$.
\item We say that $\mathcal{N}$ is a \emph{blunt rank extension} of $\mathcal{M}$, and write $\mathcal{M} \subseteq_{\mathsf{blunt}}^{\mathrm{rk}} \mathcal{N}$, if $\mathcal{M} \subseteq_e^{\mathrm{rk}} \mathcal{N}$, $M \neq N$ and $\mathcal{N}$ is not a topless rank extension of $\mathcal{M}$.
\end{itemize}

Note that $\mathsf{KP}^-$ is a subtheory of $\mathsf{M}^-+\Delta_0\textrm{-\textsf{Collection}}$. We will make use of the following results:
\begin{itemize}
\item A consequence of \cite[Theorem Scheme 6.9(i)]{mat01} is that $\mathsf{M}$ proves $\Delta_0^\mathcal{P}$-\textsf{Separation}.
\item The availability of the collection scheme for the relevant class of formulae means that the class of formulae that are equivalent to a $\Sigma_1$-formula and the class of formulae that are equivalent to a $\Pi_1$-formula are closed under bounded quantification in the theory $\mathsf{KP}^-$; the class of formulae equivalent to a $\Sigma_1^\mathcal{P}$-formula and the class of formulae that are equivalent to a $\Pi_1^\mathcal{P}$-formula are closed under bounded quantification in the theory $\mathsf{M}^-+\Delta_0^\mathcal{P}\textrm{-\textsf{Collection}}$; the class of formulae equivalent to a $\Sigma_1(\mathcal{L}^*)$-formula and the class of formulae that are equivalent to a $\Pi_1(\mathcal{L}^*)$-formula are closed under bounded quantification in the theory $\mathsf{KPU}_{\mathbb{C}\mathrm{ov}}$; and the class of formulae equivalent to a $\Sigma_1(\mathcal{L}_{\mathsf{P}}^{*})$-formula and the class of formulae that are equivalent to a $\Pi_1(\mathcal{L}_{\mathsf{P}}^{*})$-formula are closed under bounded quantification in the theory $\mathsf{KPU}^{\mathcal{P}}_{\mathbb{C}\mathrm{ov}}$.
\item The proof of \cite[I.4.4]{bar75} shows:
\begin{enumerate}
\item $\mathsf{KP}^{-} \vdash\Sigma_1\text{-}\mathsf{Collection}$; \item $\mathsf{M}^{-}+\Delta_0^{\mathcal{P}}\text{-}\mathsf{Collection} \vdash  \Sigma_1^{\mathcal{P}}\text{-}\mathsf{Collection}$; \item $\mathsf{KPU}_{\mathbb{C}\mathrm{ov}} \vdash \Sigma_1(\mathcal{L}^*)\text{-}\mathsf{Collection}$; and \item $\mathsf{KPU}^{\mathcal{P}}_{\mathbb{C}\mathrm{ov}} \vdash \Sigma_1(\mathcal{L}_{\mathsf{P}}^{*})\text{-}\mathsf{Collection}$.
\end{enumerate}
\item The argument used in \cite[I.4.5]{bar75} shows:
\begin{enumerate}
\item $\mathsf{KP}^{-} \vdash \Delta_1\text{-}\mathsf{Separation}$;

\item $\mathsf{M}^{-}+\Delta_{0}^\mathcal{P}\textrm{-}\mathsf{Collection}\vdash \Delta_1^\mathcal{P}\text{-}\mathsf{Separation}$; \item $\mathsf{KPU}_{\mathbb{C}\mathrm{ov}} \vdash\Delta_1(\mathcal{L}^*)\text{-}\mathsf{Separation}$; and \item $\mathsf{KPU}^{\mathcal{P}}_{\mathbb{C}\mathrm{ov}} \vdash \Delta_1(\mathcal{L}_{\mathsf{P}}^{*})\text{-}\mathsf{Separation}$.
\end{enumerate}
\end{itemize}

The following is Mathias's calibration \cite[Proposition Scheme 6.12]{mat01} of \cite[Theorem 6]{tak72}.

\begin{Theorems1} \label{Th:TakahashiRelationships} The following inclusions hold between the indicated classes of formulae $( n \geq 1)$:

\begin{enumerate}
\item $\Sigma_1 \subseteq (\Delta_1^\mathcal{P})^{\mathsf{MOST}}$ and $\Delta_0^\mathcal{P} \subseteq \Delta_2^{\mathbf{S}_1}$.

\item $\Sigma_{n+1} \subseteq (\Sigma_n^\mathcal{P})^{\mathsf{MOST}}.$

\item $\Pi_{n+1} \subseteq (\Pi_n^\mathcal{P})^{\mathsf{MOST}}.$

\item $\Delta_{n+1} \subseteq (\Delta_n^\mathcal{P})^{\mathsf{MOST}}.$ \item $\Sigma_n^\mathcal{P} \subseteq \Sigma_{n+1}^{\mathbf{S}_1}.$
\item $\Pi_n^\mathcal{P} \subseteq \Pi_{n+1}^{\mathbf{S}_1}.$

\item $\Delta_n^\mathcal{P} \subseteq \Delta_{n+1}^{\mathbf{S}_1}$.
\end{enumerate}
\end{Theorems1}
\begin{itemize}

\item As noted by Mathias in \cite[Corollary 6.15]{mat01}, $\mathsf{MOST}+\Pi_1\textrm{-\textsf{Collection}}$ and $\mathsf{MOST}+\Delta_0^\mathcal{P}\textrm{-\textsf{Collection}}$ axiomatise the same theory. This fact follows from part 1 of Theorem \ref{Th:TakahashiRelationships} and the results mentioned above and will be repeatedly used throughout this paper.

\end{itemize}
The $\Sigma_1^\mathcal{P}$-Recursion Theorem \cite[Theorem 6.26]{mat01} shows that the theory $\mathsf{KP}^\mathcal{P}$ is capable of constructing the levels of the cumulative hierarchy:
$$V_0= \emptyset \textrm{ and for all ordinals } \alpha,$$
$$V_{\alpha+1}= \mathcal{P}(V_\alpha) \textrm{ and, if } \alpha \textrm{ is a limit ordinal, } V_\alpha= \bigcup_{\beta \in \alpha} V_\beta.$$
More precisely, let $\mathsf{RK}(\alpha, f)$ be the $\mathcal{L}$-formula:
$$\begin{array}{c}
(f \textrm{ is a function}) \land (\alpha \textrm{ is an ordinal}) \land \mathrm{dom}(f)= \alpha ~\land\\
(\forall \beta \in \alpha)\left(\begin{array}{c}
\left((\beta \textrm{ is a limit ordinal}) \Rightarrow f(\beta)= \bigcup_{\gamma \in \beta} f(\gamma)\right) \land\\
(\exists \gamma \in \beta)(\beta= \gamma+1) \Rightarrow ((\forall x \subseteq f(\gamma))(x \in f(\beta)) \land (\forall x \in f(\beta))(x \subseteq f(\gamma)))
\end{array}\right).
\end{array}$$
Note that $\mathsf{RK}(f, \alpha)$ is a $\Delta_0^\mathcal{P}$-formula.

\begin{Lemma1}
The theory $\mathsf{KP}^\mathcal{P}$ proves
\begin{itemize}
\item[(I)] for all ordinals $\alpha$, there exists $f$ such that $\mathsf{RK}(\alpha, f)$;
\item[(II)] for all ordinals $\alpha$ and for all $f$, if $\mathsf{RK}(f, \alpha+1)$, then
$$f(\alpha)= \{x \mid \rho(x) < \alpha \}.$$
\end{itemize}
\end{Lemma1}

\noindent Therefore, the theory $\mathsf{KP}^\mathcal{P}$ proves that the function $\alpha \mapsto V_\alpha$ is total and that the graph of this function is $\Delta_1^\mathcal{P}$-definable.

In contrast, the theory $\mathsf{MOST}$ does not prove that the function $\alpha \mapsto V_\alpha$ is total (see Example \ref{Ex:PowersetPreservingExtensionsThatAreNotRank1} below). Note that the availability of \textsf{AC} in $\mathsf{MOST}$ allows us to identify cardinals with initial ordinals. Consider the $\Delta_0^\mathcal{P}$-formula $\mathsf{BFEXT}(R, X)$ defined\footnote{The abbreviation \textsf{BFEXT} has long been used by NF-theorists for well-founded extensional relations with a top, it is an abbreviation of \emph{Bien Fond\'ee Extensionnelle}, extensively employed by the French-speaking NF-ists in Belgium.} by:
$$\begin{array}{c}
(R \textrm{ is an extensional relation on } X \textrm{ with a top element})~\land\\
(\forall S \subseteq X)(S \neq \emptyset \Rightarrow (\exists x \in S)(\forall y \in S)(\langle y, x \rangle \notin R)).
\end{array}$$
The following lemma captures two important features of the theory $\mathsf{MOST}$ that follow from \cite[Theorem 3.18]{mat01}.

\begin{Lemma1} \label{Th:KeyConsequencesOfMOST}
The theory $\mathsf{MOST}$ proves the following statements:
\begin{itemize}
\item[(I)] for all $\langle X, R\rangle$ with $\mathsf{BFEXT}(X, R)$, there exists a transitive set $T$ such that $\langle X, R\rangle \cong \langle T, \in \rangle$;
\item[(II)] there exist arbitrarily large initial ordinals;
\item[(III)] for all cardinals $\kappa$, the set $H_{\leq \kappa}= \{x \mid |\mathrm{TC}(x)|\leq \kappa\}$ exists.
\end{itemize}
\end{Lemma1}

In the theory $\mathsf{MOST}$, the formula ``$X=H_{\leq \kappa}$" is $\Delta_1^\mathcal{P}$ with parameters $X$ and $\kappa$:
$$\begin{array}{c}
(\kappa \textrm{ is a cardinal}) \land\\
(\forall R \subseteq \kappa \times \kappa)(\mathsf{BFEXT}(R, \kappa) \Rightarrow (\exists x, f, T \in X)(T= \mathsf{TC}(\{x\}) \land f: R \cong \in \upharpoonright T))\land\\
(\forall x \in X)(\exists T, f \in X)(T= \mathsf{TC}(\{x\}) \land (f: T \longrightarrow \kappa \textrm{ is injective}))
\end{array}.$$

The next result is a special case of \cite[Corollary 6.11]{gor20}:

\begin{Lemma1} \label{Th:PowersetPreservingEndExtensionsOfKPPAreRankExtensions}
Let $\mathcal{M}$ and $\mathcal{N}$ be models of $\mathsf{KP}^\mathcal{P}$. If $\mathcal{M} \subseteq_e^\mathcal{P} \mathcal{N}$, then $\mathcal{M} \subseteq_e^{\mathsf{rk}} \mathcal{N}$. \Square
\end{Lemma1}

The following examples show that neither of assumptions that $\mathcal{M}$ in Lemma \ref{Th:PowersetPreservingEndExtensionsOfKPPAreRankExtensions} satisfies $\Delta_0^\mathcal{P}$-\textsf{Collection} and $\Pi_1^\mathcal{P}$-\textsf{Foundation} can be removed. The structure $\mathcal{M}$ defined in Example \ref{Ex:PowersetPreservingExtensionsThatAreNotRank1} satisfies all of the axioms of $\mathsf{KP}^{\mathcal{P}}$ except $\Delta_0^\mathcal{P}$-\textsf{Collection}. The structure $\mathcal{M}$ defined in Example \ref{Ex:PowersetPreservingExtensionsThatAreNotRank2} satisfies all of the axioms of $\mathsf{KP}^{\mathcal{P}}$ except $\Pi_1^\mathcal{P}$-\textsf{Foundation}.

\begin{Examp1} \label{Ex:PowersetPreservingExtensionsThatAreNotRank1}
Let $\mathcal{N}= \langle N, \mathsf{E}^\mathcal{N} \rangle \models\mathsf{ZF+V=L}$, and $$\mathcal{M}= \langle (H_{\aleph_\omega}^\mathcal{N})^*, \mathsf{E}^\mathcal{N} \rangle.$$
Then $\mathcal{M}\models\mathsf{MOST}+\Pi_\infty\text{-}\mathsf{Separation}$, and $\mathcal{M} \subseteq_{\mathsf{blunt}}^\mathcal{P} \mathcal{N}$, but $\mathcal{N}$ is not a rank extension of $\mathcal{M}$.
\end{Examp1}

\begin{Examp1} \label{Ex:PowersetPreservingExtensionsThatAreNotRank2}
Let $\mathcal{N}= \langle N, \mathsf{E}^\mathcal{N} \rangle$ be an $\omega$-nonstandard model of $\mathsf{ZF+V=L}$. Let $\mathcal{M}= \langle M, \mathsf{E}^\mathcal{N} \rangle$, where
$$M= \bigcup_{n \in \omega} (H_{\aleph_n}^\mathcal{N})^*.$$
Then $\mathcal{M}\models\mathsf{MOST}+\Pi_1\text{-}\mathsf{Collection}$, and $\mathcal{M} \subseteq_{\mathsf{topless}}^\mathcal{P} \mathcal{N}$, but $\mathcal{N}$ is not a rank extension of $\mathcal{M}$.
\end{Examp1}

 The following recursive definition can be carried out within $\mathsf{KPU}_{\mathbb{C}\mathrm{ov}}$ thanks to the ability of  $\mathsf{KPU}_{\mathbb{C}\mathrm{ov}}$ to carry out $\Sigma_1(\mathcal{L}^*)$-recursions. The recursion defines an operation $\ulcorner \cdot \urcorner$ for coding the infinitary formulae of $\mathcal{L}^{\mathsf{ee}}_{\infty \omega}$, where  $\mathcal{L}^{\mathsf{ee}}$ be the language obtained from $\mathcal{L}$ by adding new constant symbols $\bar{a}$ for each urelement $a$ and a new constant symbol $\mathbf{c}$.
\begin{itemize}
\item for all ordinals $\alpha$, $\ulcorner v_\alpha \urcorner= \langle 0, \alpha \rangle$,


\item for all urelements $a$, $\ulcorner \bar{a} \urcorner= \langle 1, a \rangle$,
\item $\ulcorner \mathbf{c} \urcorner= \langle 2, 0 \rangle$,
\item if $\phi$ is an $\mathcal{L}^{\mathsf{ee}}_{\infty \omega}$-formula and $x$ is a free variable of $\phi$, then
$$\ulcorner \exists x \phi \urcorner= \langle 3, \ulcorner x \urcorner, \ulcorner \phi \urcorner \rangle,$$
\item if $\phi$ is an $\mathcal{L}^{\mathsf{ee}}_{\infty \omega}$-formula and $x$ is a free variable of $\phi$, then
$$\ulcorner \forall x \phi \urcorner= \langle 4, \ulcorner x \urcorner, \ulcorner \phi \urcorner \rangle,$$
\item if $\Phi$ is a set of $\mathcal{L}^{\mathsf{ee}}_{\infty \omega}$-formulae such that only finitely many variables appear as a free variable of some formula in $\Phi$, then
$$\ulcorner \bigvee_{\phi \in \Phi} \phi \urcorner= \langle 5, \Phi^*\rangle, \textrm{ where } \Phi^*= \{ \ulcorner \phi \urcorner \mid \phi \in \Phi\},$$
\item if $\Phi$ is a set of $\mathcal{L}^{\mathsf{ee}}_{\infty \omega}$-formulae such that only finitely many variables appear as a free variable of some formula in $\Phi$, then
$$\ulcorner \bigwedge_{\phi \in \Phi} \phi \urcorner= \langle 6, \Phi^*\rangle, \textrm{ where } \Phi^*= \{ \ulcorner \phi \urcorner \mid \phi \in \Phi\},$$
\item if $\phi$ is an $\mathcal{L}^{\mathsf{ee}}_{\infty \omega}$-formula, then $\ulcorner \neg \phi \urcorner= \langle 7, \ulcorner \phi \urcorner \rangle$,
\item if $s$ and $t$ are terms of $\mathcal{L}^{\mathsf{ee}}_{\infty \omega}$, then $\ulcorner s = t \urcorner= \langle 8, \ulcorner s \urcorner, \ulcorner t \urcorner \rangle$,
\item if $s$ and $t$ are terms of $\mathcal{L}^{\mathsf{ee}}_{\infty \omega}$, then $\ulcorner s \in t \urcorner= \langle 9, \ulcorner s \urcorner, \ulcorner t \urcorner \rangle$.
\end{itemize}
Let $\mathcal{M}= \langle M, \mathsf{E}^\mathcal{M} \rangle$ be an $\mathcal{L}$-structure and let $\mathfrak{A}_\mathcal{M}= \langle \mathcal{M}; A, \in, \mathsf{F}^\mathfrak{A}, \mathsf{P}^\mathfrak{A} \rangle$ be a power admissible set covering $\mathcal{M}$. We use $\mathcal{L}^{\mathsf{ee}}_{\mathfrak{A}_\mathcal{M}}$ to denote the fragment of $\mathcal{L}^{\mathsf{ee}}_{\infty \omega}$ that is coded in $\mathfrak{A}_\mathcal{M}$. The $\mathcal{L}^{\mathsf{ee}}_{\mathfrak{A}_\mathcal{M}}$-formulae in the form $s=t$ or $s \in t$, where $s$ and $t$ are $\mathcal{L}^{\mathsf{ee}}_{\mathfrak{A}_\mathcal{M}}$-terms, are the {\bf atomic formulae} of $\mathcal{L}^{\mathsf{ee}}_{\mathfrak{A}_\mathcal{M}}$. The formula that identifies the codes of the atomic formulae of $\mathcal{L}^{\mathsf{ee}}_{\mathfrak{A}_\mathcal{M}}$ is $\Delta_1(\mathcal{L}^*)$-definable over $\mathfrak{A}_\mathcal{M}$. Similarly, other important properties of codes of $\mathcal{L}^{\mathsf{ee}}_{\mathfrak{A}_\mathcal{M}}$ constituents, such as being a {\it variable}, {\it constant}, {\it well-formed formula}, {\it sentence}, \ldots, are all $\Delta_1(\mathcal{L}^*)$-definable over $\mathfrak{A}_\mathcal{M}$. We will often equate an $\mathcal{L}^{\mathsf{ee}}_{\mathfrak{A}_\mathcal{M}}$-theory $T$ with that subset of $\mathfrak{A}_\mathcal{M}$ of codes of $\mathcal{L}^{\mathsf{ee}}_{\mathfrak{A}_\mathcal{M}}$-sentences in $T$.

\medskip

The following is the Barwise Compactness Theorem (\cite[III.5.6]{bar75}) tailormade for countable admissible $\mathcal{L}_{\mathsf{P}}^{*}$-structures.

\begin{Theorems1} \label{Th:BarwiseCompactnessTheorem}
(Barwise Compactness Theorem) Let $\mathfrak{A}_\mathcal{M}= \langle \mathcal{M}; A, \in, \mathsf{F}^\mathfrak{A}, \mathsf{P}^\mathfrak{A} \rangle$ be a power admissible set covering $\mathcal{M}$. Let $T$ be an $\mathcal{L}^{\mathsf{ee}}_{\mathfrak{A}_\mathcal{M}}$-theory that is $\Sigma_1(\mathcal{L}_{\mathsf{P}}^{*})$-definable over $\mathfrak{A}_\mathcal{M}$ and such that for all $T_0 \subseteq T$, if $T_0 \in A$, then $T_0$ has a model. Then $T$ has a model.
\end{Theorems1}

\section[Sigma 1 in powerset \textsf{Foundation}]{The scheme of $\Sigma_1^\mathcal{P}$-Foundation}

Motivated by the apparent reliance of the constructions presented in the next section on $\Sigma_1^\mathcal{P}$-\textsf{Foundation}, this section investigates the status of this scheme in the theories $\mathsf{MOST}+\Pi_1\textrm{-\textsf{Collection}}$ and $\mathsf{KP}^\mathcal{P}$. We begin by showing that $\mathsf{KP}^\mathcal{P}$ proves $\Sigma_1^\mathcal{P}$-\textsf{Foundation}. In contrast, $\Sigma_1^\mathcal{P}$-\textsf{Foundation} is not provable in $\mathsf{MOST}+\Pi_1\textrm{-\textsf{Collection}}$ but does hold in every $\omega$-standard model of this theory.

In \cite[Lemma 4.4]{rat20} it is shown that $\mathsf{KP}^\mathcal{P}+\mathsf{AC}$ proves $\Sigma_1^\mathcal{P}$-\textsf{Foundation}\footnote{Rathjen proves the scheme that asserts that set induction holds for all $\Pi_1^\mathcal{P}$-formulae, which, in the theory $\mathsf{KP}^\mathcal{P}$, is equivalent to $\Sigma_1^\mathcal{P}$-\textsf{Foundation}.}. Here we use a modification of a choiceless scheme of dependant choices introduced in \cite{rat92} to show that $\Sigma_1^\mathcal{P}$-\textsf{Foundation} can be proved in $\mathsf{KP}^\mathcal{P}$. The following is \cite[Definition 3.1]{rat92}:

\begin{Definitions1}
Let $\phi(x, y, \vec{z})$ be an $\mathcal{L}$-formula. Define $\delta^\phi(a, b, f, \vec{z})$ to be the formula:
$$a \textrm{ is an ordinal} \Rightarrow \left( \begin{array}{c}
f \textrm{ is a function} \land \mathsf{dom}(f)= a+1 \land f(0)= \{b\}\land\\
(\forall u \in a)\left( \begin{array}{c}
(\forall x \in f(u))(\exists y \in f(u+1)) \phi(x, y, \vec{z}) \land\\
(\forall y \in f(u+1))(\forall x \in f(u)) \phi(x, y, \vec{z}) \end{array}\right)
\end{array}\right).$$
\end{Definitions1}

By considering the variables $\vec{z}$ to be parameters, $\phi(x, y, \vec{z})$ defines a directed graph. The formula $\delta^\phi(a, b, f, \vec{z})$ says that for all $0 \leq i \leq a$, $f(i)$ is a collection of vertices lying at a stage $i$ on a directed path of length $a$ starting at $b$ in this graph. In the next definition we introduce a formula that, given $b$, $f$ and $\vec{z}$, says that $f$ is a function with domain $\omega$ and for all $n \in \omega$, $\delta^\phi(n+1, b, f \upharpoonright (n+1), \vec{z})$.

\begin{Definitions1}
Let $\phi(x, y, \vec{z})$ be an $\mathcal{L}$-formula. Define $\delta_\omega^\phi(b, f, \vec{z})$ to be the formula:
$$\begin{array}{c}
f \textrm{ is a function} \land \mathsf{dom}(f)= \omega \land f(0)= \{b\}\land\\
(\forall u \in \omega)\left( \begin{array}{c}
(\forall x \in f(u))(\exists y \in f(u+1)) \phi(x, y, \vec{z}) \land\\
(\forall y \in f(u+1))(\forall x \in f(u)) \phi(x, y, \vec{z}) \end{array}\right)
\end{array}.$$
\end{Definitions1}

Note that if $\phi(x, y, \vec{z})$ is a $\Delta_0^\mathcal{P}$-formula ($\Delta_0$-formula) then both $\delta^\phi(a, b, f, \vec{z})$ and $\delta_\omega^\phi(b, f, \vec{z})$ are both $\Delta_0^\mathcal{P}$-formulae (respectively $\Delta_0$-formulae). The following is a modification of Rathjen's $\Delta_0$-weak dependant choices scheme ($\Delta_0\text{-}\mathsf{WDC}$) from \cite{rat92}:
\begin{itemize}
\item[]($\Delta_0^\mathcal{P}\text{-}\mathsf{WDC}_\omega$) For all $\Delta_0^\mathcal{P}$-formulae, $\phi(x, y, \vec{z})$,
$$\forall \vec{z}(\forall x \exists y \phi(x, y, \vec{z}) \Rightarrow \forall w \exists f \delta_\omega^\phi(w, f, \vec{z})).$$
\end{itemize}

The next result is based on the proof of \cite[Proposition 3.2]{rat92}:

\begin{Theorems1} \label{Th:Delta0PWDCProvesSigma1PFoundation}
The theory $\mathsf{M}+\Delta_0^\mathcal{P}\text{-}\mathsf{WDC}_\omega$ proves $\Sigma_1^\mathcal{P}\text{-}\mathsf{Foundation}$.
\end{Theorems1}

\begin{proof}
Work in the theory $\mathsf{M}+\Delta_0^\mathcal{P}\text{-}\mathsf{WDC}_\omega$. Suppose, for a contradiction, that there is an instance of $\Sigma_1^\mathcal{P}$-\textsf{Foundation} that fails. Let $\phi(x, y, \vec{z})$ be a $\Delta_0^\mathcal{P}$-formula and let $\vec{a}$ be a finite sequence of sets such that the class $C= \{x \mid \exists y \phi(x, y, \vec{a})\}$ is nonempty and has no $\in$-least element. Let $b$ and $d$ be such that $\phi(b, d, \vec{a})$ holds. Now, since $C$ has no $\in$-least element,
$$\forall x \forall u \exists y \exists v(\phi(x, u, \vec{a}) \Rightarrow (y \in x) \land \phi(y, v, \vec{a})).$$
Therefore, we have have $\forall x \exists y \theta(x, y, \vec{a})$ where $\theta(x, y, \vec{a})$ is
$$x= \langle x_0, x_1 \rangle \land y= \langle y_0, y_1 \rangle \land (\phi(x_0, x_1, \vec{a}) \Rightarrow (y_0 \in x_0) \land \phi(y_0, y_1, \vec{a})).$$
Note that $\theta(x, y, \vec{a})$ is a $\Delta_0^\mathcal{P}$-formula. Therefore, using $\Delta_0^\mathcal{P}\text{-}\mathsf{WDC}_\omega$, let $f$ be such that $\delta_\omega^\theta(\langle b, d \rangle, f, \vec{a})$. Now, $\Delta_0^\mathcal{P}$-\textsf{Separation} facilitates induction for $\Delta_0^\mathcal{P}$-formulae and proves that for all $n \in \omega$,
$$\begin{array}{c}
f(n)\neq \emptyset \land (\forall x \in f(n))(x= \langle x_0, x_1 \rangle \land \phi(x_0, x_1, \vec{a})) \land\\
(\forall x \in f(n))(\exists y \in f(n+1))(x= \langle x_0, x_1 \rangle \land y= \langle y_0, y_1 \rangle \land y_0 \in x_0)\land\\
(\forall y \in f(n+1))(\exists x \in f(n))(y= \langle y_0, y_1\rangle \land x= \langle x_0, x_1 \rangle \land y_0 \in x_0)
\end{array}.$$
Let $B= \mathsf{TC}(\{b\})$. Induction for $\Delta_0$-formulae suffices to prove that for all $n \in \omega$,
$$(\forall x \in f(n))(x= \langle x_0, x_1 \rangle \land x_0 \in B).$$
Consider
$$A= \left\{ x \in B \mid (\exists n \in \omega)(\exists z \in f(n)) \left(\exists y \in \bigcup z\right)(z= \langle x, y \rangle) \right\},$$
which is a set by $\Delta_0$-\textsf{Separation}. Now, let $x \in A$. Let $y$ and $n \in \omega$ be such that $\langle x, y \rangle \in f(n)$. Therefore, there exists $w \in f(n+1)$ such that $w= \langle u, v \rangle$ and $u \in x$. So $u \in A$ and $u \in x$, which shows that $A$ has no $\in$-least element. This contradicts \textsf{Set-Foundation} in $\mathsf{M}$ and proves the theorem.
\Square
\end{proof}

The fact that $\mathsf{KP}^\mathcal{P}$ proves $\Sigma_1^\mathcal{P}$-\textsf{Foundation} follows from the fact that $\mathsf{KP}^\mathcal{P}$ proves $\Delta_0^\mathcal{P}\text{-}\mathsf{WDC}_\omega$. The proof of Theorem \ref{Th:KPPProvesDelta0PWDC} is inspired by the argument used in the proof of \cite[Theorem 4.15]{flw16}. The stratification of the universe into ranks allows us to select sets of paths through a relation defined by a $\Delta_0^\mathcal{P}$-formula $\phi(x, y, \vec{z})$ with parameters $\vec{z}$.

\begin{Definitions1}
Let $\phi(x, y, \vec{z})$ be an $\mathcal{L}$-formula. Define $\eta^\phi(a, b, f, \vec{z})$ by
$$\begin{array}{c}
\delta^\phi(a, b, f, \vec{z}) \land\\
(\forall u \in a) \exists \alpha \exists X \left(\begin{array}{c}
(\alpha \textrm{ is an ordinal}) \land (X= V_\alpha) ~ \land\\
(\forall x \in f(u+1))(x \in X)~ \land\\
(\forall y \in X)(\forall x \in f(u))(\phi(x, y, \vec{z}) \Rightarrow y \in f(u+1))~ \land\\
(\forall \beta \in \alpha)(\forall Y \in X)\left(\begin{array}{c}
Y= V_\beta \Rightarrow\\
(\exists x \in f(u))(\forall y \in Y) \neg \phi(x, y, \vec{z})
\end{array}\right)
\end{array}\right)
\end{array}.$$
\end{Definitions1}

The formula $\eta^\phi(a, b, f, \vec{z})$ asserts that $f$ is a function with domain $a+1$ such that $f(0)= \{b\}$ and for all $u \in a$, $f(u+1)$ is the set of $y$ of rank $\alpha$ such that there exists $x \in f(u)$ with $\phi(x, y, \vec{z})$ and $\alpha$ is the minimal ordinal such that for all $x \in f(u)$, there exists $y$ of rank $\alpha$ such that $\phi(x, y, \vec{z})$. Recall that, in the theory $\mathsf{KP}^\mathcal{P}$, the formula `$X=V_\alpha$' is $\Delta_1^\mathcal{P}$ with parameters $X$ and $\alpha$. Therefore, if $\phi(x, y, \vec{z})$ is a $\Delta_0^\mathcal{P}$-formula, then $\eta^\phi(a, b, f, \vec{z})$ is equivalent to a $\Sigma_1^\mathcal{P}$-formula in the theory $\mathsf{KP}^\mathcal{P}$.

\begin{Theorems1} \label{Th:KPPProvesDelta0PWDC}
The theory $\mathsf{KP}^\mathcal{P}$ proves $\Delta_0^\mathcal{P}\text{-}\mathsf{WDC}_\omega$.
\end{Theorems1}

\begin{proof}
Work in the theory $\mathsf{KP}^\mathcal{P}$. Let $\phi(x, y, \vec{z})$ be a $\Delta_0^\mathcal{P}$-formula. Let $\vec{a}$ be sets such that $\forall x \exists y \phi(x, y, \vec{a})$ holds. Let $b$ be a set. We begin by claiming that for all $n \in \omega$, $\exists f \eta^\phi(n, b, f, \vec{a})$. Suppose, for a contradiction, that this does not hold. Using $\Pi_1^\mathcal{P}$-\textsf{Foundation}, there exists a least $m \in \omega$ such that $\neg \exists f \eta^\phi(m, b, f, \vec{a})$. It is straightforward to see that $m \neq 0$. Therefore, there exists a function $g$ with $\mathsf{dom}(g)= m$ such that $\eta^\phi(m-1, b, g, \vec{a})$ holds. Consider the class
$$A= \{\alpha \in \mathsf{Ord} \mid \forall X (X= V_\alpha \Rightarrow (\forall x \in g(m-1))(\exists y \in X)\phi(x, y, \vec{a}))\}.$$
Applying $\Delta_0^\mathcal{P}$-\textsf{Collection} to the formula $\phi(x, y, \vec{a})$ shows that $A$ is nonempty. Therefore, by $\Pi_1^\mathcal{P}$-\textsf{Foundation}, there exists a least element $\beta \in A$. Let
$$C= \{ y \in V_\beta \mid (\exists x \in g(m-1))\phi(x, y, \vec{a})\},$$
which is a set by $\Delta_0^\mathcal{P}$-\textsf{Separation}. Now, let $f= g \cup \{\langle m, C \rangle\}$. So, $\eta^\phi(m, b, f, \vec{a})$, which is a contradiction. Therefore, for all $n \in \omega$, $\exists f \eta^\phi(n, b, f, \vec{a})$. Note that for all $n \in \omega$ and for all $f$ and $g$, if $\eta^\phi(n, b, f, \vec{a})$ and $\eta^\phi(n, b, g, \vec{a})$, then $f=g$. Now, using $\Sigma_1^\mathcal{P}$-\textsf{Collection}, we can find a set $D$ such that $(\forall n \in \omega)(\exists f \in D)\eta^\phi(n, b, f, \vec{a})$. Let
$$f= \{ \langle n, X \rangle \in \omega \times \mathsf{TC}(D) \mid (\exists g \in D)(\eta^\phi(n, b, g, \vec{a})\land g(n)=X) \}$$
$$=\{ \langle n, X \rangle \in \omega \times \mathsf{TC}(D) \mid (\forall g \in D)(\eta^\phi(n, b, g, \vec{a}) \Rightarrow g(n)= X) \}.$$
Now, $f$ is a set by $\Delta_1^\mathcal{P}$-Separation and $f$ is the function required by $\Delta_0^\mathcal{P}\text{-}\mathsf{WDC}_\omega$.
\Square
\end{proof}

Combining Theorems \ref{Th:Delta0PWDCProvesSigma1PFoundation} and \ref{Th:KPPProvesDelta0PWDC} yields:

\begin{Coroll1} \label{Th:KPPProvesSigma1PFoundation}
$\mathsf{KP}^\mathcal{P} \vdash \Sigma_1^\mathcal{P}\text{-}\mathsf{Foundation}$. \Square
\end{Coroll1}

We now turn to investigating $\Sigma_1^\mathcal{P}$-\textsf{Foundation} in the theory $\mathsf{MOST}+\Pi_1\textrm{-\textsf{Collection}}$. The following is an instance of \cite[Proposition 2]{pk78} in the context of set theory:

\begin{Theorems1} \label{Th:ParisKirbyResult}

Let $\Sigma$ denote $\Sigma_1^\mathcal{P}\text{-}\mathsf{Foundation~on~}  \omega$, and $\Pi$ denote $\Pi_1^\mathcal{P}\text{-}\mathsf{Foundation~on~}  \omega$. Then we have:
\begin{center}
$\mathsf{M}^-+\Delta_0^\mathcal{P}\text{-}\mathsf{Collection} + \Pi\vdash \Sigma$, and $\mathsf{M}^- + \Sigma \vdash \Pi$.
\end{center}


\end{Theorems1}

\begin{proof}
To see that $\Pi$ implies $\Sigma$, work in the theory $\mathsf{M}^-+\Delta_0^\mathcal{P}\textrm{-}\mathsf{Collection}$. We prove the contrapositive. Let $\phi(x, \vec{z})$ be a $\Pi_1^\mathcal{P}$-formula and let $\vec{a}$ be sets such that the class $\{x \in \omega \mid \phi(x, \vec{a})\}$ is nonempty and has no least element. Let $p \in \omega$ be such that $\phi(p, \vec{a})$. Let
$$C= \{ x \in \omega \mid \exists w (x+w=p \land (\forall y \in w)\neg \phi(y, \vec{a}))\}.$$
Note that $\Delta_0^\mathcal{P}$-\textsf{Collection} implies that $C$ is a $\Sigma_1^\mathcal{P}$-definable subclass of $\omega$. Moreover, $p \in C$ and $0 \notin C$. Identical reasoning to that used above shows that $C$ has no least element. Therefore $\Sigma_1^\mathcal{P}\text{-}\mathsf{Foundation}$ on $\omega$ fails.

To see that $\Sigma$ implies $\Pi$, work in the theory $\mathsf{M}^-$. Again, we prove the contrapositive. Let $\phi(x, \vec{z})$ be a $\Sigma_1^\mathcal{P}$-formula and let $\vec{a}$ be the sequence of set parameters such that the class $\{x \in \omega \mid \phi(x, \vec{a})\}$ is nonempty and has no least element. Let $p \in \omega$ be such that $\phi(p, \vec{a})$. Let
$$C= \{ x \in \omega \mid \forall w(x+w= p \Rightarrow (\forall y \in w)\neg \phi(y, \vec{a})\}.$$
Note that $C$ is a $\Pi_1^\mathcal{P}$-definable subclass of $\omega$, $p \in C$ and $0 \notin C$. Suppose that $q \in C$ is a least element of $C$. Let $u \in \omega$ be such that $q+u= p$. Now, $\phi(u, \vec{a})$, since $q$ is the least of $C$, and $(\forall y \in u)\neg \phi(y, \vec{a})$. But then $u$ is a least element of $\{x \in \omega \mid \phi(x, \vec{a})\}$, which is a contradiction. Therefore $\Pi_1^\mathcal{P}\text{-}\mathsf{Foundation~on~ }\omega$ fails. \Square
\end{proof}

An examination of the proof of \cite[Proposition 9.22]{mat01} yields:

\begin{Theorems1}
The consistency of $\mathsf{Mac}$ is provable in $\mathsf{M}+\Pi_1^\mathcal{P}\text{-}\mathsf{Foundation~ on~}\omega$. \Square
\end{Theorems1}

The results of \cite{mat01} and \cite{mck19} (see \cite[Corollary 3.5]{mck19}) show that $\mathsf{Mac}$ and $\mathsf{MOST}+\Pi_1\text{-}\mathsf{Collection}$ have the same consistency strength. Therefore, Theorem \ref{Th:ParisKirbyResult} yields:

\begin{Theorems1}
The consistency of $\mathsf{MOST}+\Pi_1\text{-}\mathsf{Collection}$ is provable in $\mathsf{MOST}+\Pi_{1}\text{-}\mathsf{Collection}+\Sigma_1^{\mathcal{P}}\text{-}\mathsf{Foundation}$. \Square
\end{Theorems1}

Therefore, $\Sigma_1^\mathcal{P}$-\textsf{Foundation} is not provable in $\mathsf{MOST}+\Pi_1\text{-}\mathsf{Collection}$. However, we can show that $\Sigma_1^\mathcal{P}$-\textsf{Foundation} does hold in every model of $\mathsf{MOST}+\Pi_1\text{-}\mathsf{Collection}$ in which the natural numbers are standard.

In the context of the theory $\mathsf{MOST}+\Pi_1\textrm{-\textsf{Collection}}$, we can use the stratification of the universe into the sets $H_{\leq \kappa}$ in the same way that we used the stratification of the universe into ranks in the proof of Theorem \ref{Th:KPPProvesDelta0PWDC}.

\begin{Definitions1}
Let $\phi(x, y, \vec{z})$ be an $\mathcal{L}$-formula. Define $\chi^\phi(a, b, f, \vec{z})$ by
{\footnotesize $$\begin{array}{c}
\delta^\phi(a, b, f, \vec{z})\land\\
(\forall u \in a)\exists \kappa \exists X \left(\begin{array}{c}
(X= H_{\leq \kappa}) \land (\forall x \in f(u+1))(x \in X) \land\\
(\forall y \in X)(\forall x \in f(u))(\phi(x, y, \vec{z}) \Rightarrow y \in f(u+1)) \land\\
(\forall \lambda \in \kappa)(\exists x \in f(u))(\forall R \subseteq \lambda \times \lambda)\left(\begin{array}{c}
\mathsf{BFEXT}(R, \lambda) \Rightarrow\\
\exists T, f, y \left(\begin{array}{c}
T=\mathsf{TC}(\{y\}) \land\\
f: R \cong \in \upharpoonright T \land \neg \phi(x, y, \vec{z})
\end{array}\right)
\end{array}\right)
\end{array}\right)
\end{array}.$$}
\end{Definitions1}

\noindent The formula $\chi^\phi(a, b, f, \vec{z})$ asserts that $f$ is a function with domain $a+1$ such that $f(0)= \{b\}$ and for all $u \in a$, $f(u+1)$ is the set of all $y$ in $H_{\leq \kappa}$ such that there exists an $x \in f(u)$ with $\phi(x, y, \vec{z})$ and $\kappa$ is the minimal cardinal such that for all $x \in f(u)$, there exists $y \in H_{\leq \kappa}$ with $\phi(x, y, \vec{z})$. Recall that the formula expressing ``$X=H_{\leq \kappa}$" is $\Delta_1^\mathcal{P}$ with parameters $X$ and $\kappa$ in the theory $\mathsf{MOST}$. Therefore, if $\phi(x, y, \vec{z})$ is a $\Delta_0^\mathcal{P}$-formula, then $\chi^\phi(a, b, f, \vec{z})$ is equivalent to a $\Sigma_1^\mathcal{P}$-formula in the theory $\mathsf{MOST}+\Pi_1\textrm{-\textsf{Collection}}$. In the proof of the next theorem the formula $\chi^\phi(a, b, f, \vec{z})$ plays the role of $\eta^\phi(a, b, f, \vec{z})$ in the proof of Theorem \ref{Th:KPPProvesDelta0PWDC}.

\begin{Theorems1} \label{Th:Delta0PWDCHoldsInOmegaModels}
Let $\mathcal{M}= \langle M, \mathsf{E}^\mathcal{M} \rangle$ be an $\omega$-standard model of $\mathsf{MOST}+\Pi_1\text{-}\mathsf{Collection}$. Then $\mathcal{M} \models \Delta_0^\mathcal{P}\text{-}\mathsf{WDC}_\omega$.
\end{Theorems1}

\begin{proof}
Let $\phi(x, y, \vec{z})$ be a $\Delta_0^\mathcal{P}$-formula. Let $\vec{a}$ be sets such that $\mathcal{M} \models \forall x \exists y \phi(x, y, \vec{a})$. Let $b$ a set. We begin by showing that
$$\mathcal{M} \models (\forall n \in \omega)\exists f \chi^\phi(n, b, f, \vec{a}).$$
Work inside $\mathcal{M}$. Suppose, for a contradiction, that there exists $n \in \omega$ such that $\neg \exists f \chi^\phi(n, b, f, \vec{a})$ holds. Therefore, since $\mathcal{M}$ is $\omega$-standard, there is a least $m \in \omega$ such that $\neg \exists f \chi^\phi(m, b, f, \vec{a})$. It is straightforward to see that $m \neq 0$. Therefore, there exists a function $g$ with $\mathsf{dom}(g)= m$ such that $\chi^\phi(m-1, b, g, \vec{a})$. Consider
$$A= \{ \kappa \in \mathsf{Ord}^\mathcal{M} \mid \mathcal{M} \models \forall X(X= H_{\leq \kappa} \Rightarrow (\forall x \in g(m-1))(\exists y \in X)\phi(x, y, \vec{a}))\}$$
$$= \{ \kappa \in \mathsf{Ord}^\mathcal{M} \mid \mathcal{M} \models \exists X(X= H_{\leq \kappa} \land (\forall x \in g(m-1))(\exists y \in X)\phi(x, y, \vec{a}))\}.$$
Applying $\Delta_0^\mathcal{P}$-\textsf{Collection} to $\phi(x, y, \vec{a})$ shows that $A$ is nonempty. Therefore, $\Delta_1^\mathcal{P}$-\textsf{Separation} ensures that $A$ has an $\in$-least element $\lambda$. Let
$$C= \{ y \in H_{\leq \lambda} \mid (\exists x \in g(m-1))\phi(x, y, \vec{a})\},$$
which is a set by $\Delta_0^\mathcal{P}$-\textsf{Separation}. Now, let $f= g \cup \{\langle m, C \rangle\}$. So, $\chi^\phi(m, b, f, \vec{a})$ which is a contradiction. This shows that
$$\mathcal{M} \models (\forall n \in \omega)\exists f \chi^\phi(n, b, f, \vec{a}).$$
Work inside $\mathcal{M}$. Note that for all $n \in \omega$ and for all $f$ and $g$, if $\chi^\phi(n, b, f, \vec{a})$ and $\chi^\phi(n, b, g, \vec{a})$, then $f=g$. Now, using $\Sigma_1^\mathcal{P}$-\textsf{Collection}, there exists $D$ such that $$(\forall n \in \omega)(\exists f \in D) \chi^\phi(n, b, f, \vec{a}).$$Let
$$f= \{\langle n, X \rangle \in \omega \times \mathsf{TC}(D) \mid (\exists g \in D)(\chi^\phi(n, b, g, \vec{a}) \land g(n)= X)\}$$
$$= \{\langle n, X \rangle \in \omega \times \mathsf{TC}(D) \mid (\forall g \in D)(\chi^\phi(n, b, g, \vec{a}) \Rightarrow g(n)= X)\},$$
which is a set by $\Delta_1^\mathcal{P}$-\textsf{Separation}. Therefore $\delta_\omega(b, f, \vec{a})$ holds, which completes the proof that $\Delta_0^\mathcal{P}\text{-}\mathsf{WDC}_\omega$ holds in $\mathcal{M}$.
\Square
\end{proof}

Combining Theorems \ref{Th:Delta0PWDCProvesSigma1PFoundation} and \ref{Th:Delta0PWDCHoldsInOmegaModels}:

\begin{Coroll1} \label{Th:Sigma1PFoundationInOmegaModelsMOST}
If $\mathcal{M}$ is an $\omega$-standard model of $\mathsf{MOST}+\Pi_1\text{-}\mathsf{Collection}$, then $\mathcal{M} \models \Sigma_1^\mathcal{P}\text{-}\mathsf{Foundation}$.\Square
\end{Coroll1}

\section[Obtaining the cover]{Obtaining $\mathbb{C}\mathsf{ov}_\mathcal{M}$ from $\mathcal{M}$} \label{Sec:Cover}

\cite[Appendix]{bar75} shows how the admissible cover, $\mathbb{C}\mathsf{ov}_\mathcal{M}$, can be built from an $\mathcal{L}$-structure $\mathcal{M}$ that satisfies $\mathsf{KP}+\Sigma_1\text{-}\mathsf{Foundation}$. The construction proceeds in two stages. The first stage interprets a model of $\mathsf{KPU}_{\mathbb{C}\mathsf{ov}}$ inside $\mathcal{M}$. The second stage takes the well-founded part of this interpreted model of $\mathsf{KPU}_{\mathbb{C}\mathsf{ov}}$ to obtain an admissible set covering $\mathcal{M}$ that \cite[Appendix]{bar75} shows is minimal. It should be noted that \cite[Appendix]{bar75} starts with a structure $\mathcal{M}$ that satisfies full $\Pi_\infty$-\textsf{Foundation}.

It is noted in \cite[Chapter 2]{res87} that all of the elements of Barwise's construction of $\mathbb{C}\mathsf{ov}_\mathcal{M}$ can be carried out when $\mathcal{M}$ satisfies $\Pi_1 \cup \Sigma_1\text{-}\mathsf{Foundation}$. The aim of this section is to review the construction of $\mathbb{C}\mathsf{ov}_\mathcal{M}$ from $\mathcal{M}$ and investigate the influence of the theory of $\mathcal{M}$ on $\mathbb{C}\mathsf{ov}_\mathcal{M}$. In particular, we will show that if $\mathcal{M}$ is a model of $\mathsf{KP}+\mathsf{powerset}+\Delta_0^\mathcal{P}\textrm{-\textsf{Collection}}+\Sigma_1^\mathcal{P}\text{-}\mathsf{Foundation}$, then $\mathsf{P}$ can be interpreted in $\mathbb{C}\mathsf{ov}_\mathcal{M}$ to make it a power admissible set.

Throughout this section we will work with a fixed $\mathcal{L}$-structure $\mathcal{M}= \langle M, \mathsf{E}^\mathcal{M} \rangle$ that satisfies $\mathsf{KP}+\mathsf{Powerset}+\Delta_0^\mathcal{P}\text{-}\mathsf{Collection}+\Sigma_1^\mathcal{P}\text{-}\mathsf{Foundation}$. We begin by expanding the interpretation of the theory $\mathsf{KPU}_{\mathbb{C}\mathsf{ov}}$ inside $\mathcal{M}$ presented in \cite[Appendix Section 3]{bar75} to obtain $\mathcal{L}_{\mathsf{P}}^{*}$-structure that satisfies \textsf{Powerset}. Working inside $\mathcal{M}$, define the unary relations $\mathsf{N}$ and $\mathsf{Set}$, the binary relations $\mathcal{E}$ and $\mathsf{E}^\prime$, and unary function symbols $\bar{\mathsf{F}}$ and $\bar{\mathsf{P}}$ by:
$$\mathsf{N}(x) \textrm{~~iff~~} \exists y (x= \langle 0, y \rangle);$$
$$x \mathsf{E}^\prime y \textrm{~~iff~~} \exists w \exists z(x = \langle 0, w \rangle \land y= \langle 0, z \rangle \land w \in z);$$
$$\mathsf{Set}(x) \textrm{~~iff~~} \exists y (x= \langle 1, y \rangle \land (\forall z \in y)(\mathsf{N}(z) \lor \mathsf{Set}(z)));$$
$$x \mathcal{E} y \textrm{~~iff~~} \exists z (y= \langle 1, z \rangle \land x \in z);$$
$$\bar{\mathsf{F}}(x)= \langle 1, X \rangle \textrm{ where } X= \{\langle 0, y \rangle \mid \exists w(x = \langle 0, w \rangle \land y \in w )\};$$
$$\bar{\mathsf{P}}(x)= \langle 1, X \rangle \textrm{ where } X= \{\langle 1, y \rangle \mid \exists w(x= \langle 1, w \rangle \land y \subseteq w)\}.$$
\cite[Appendix Section 3]{bar75} notes that $\mathsf{N}$, $\mathsf{E}^\prime$, $\mathcal{E}$ and $\bar{\mathsf{F}}$ are defined by $\Delta_0$-formulae in $\mathcal{M}$, and, using the Second Recursion Theorem (\cite[V.2.3.]{bar75}), $\mathsf{Set}$ can be expressed using a $\Sigma_1$-formula. \cite[Chapter 2]{res87} notes that the Second Recursion Theorem can be proved in $\mathsf{KP}+\Sigma_1\text{-}\mathsf{Foundation}$. The function $y= \bar{\mathsf{P}}(x)$ is defined by a $\Delta_0^\mathcal{P}$-formula:
$$y= \bar{\mathsf{P}}(x) \textrm{~~iff~~}$$
$$\mathsf{OP}(x)\land \mathsf{OP}(y)\land \mathsf{fst}(x)=1 \land \mathsf{fst}(y)= 1 \land$$
$$(\forall z \subseteq \mathsf{snd}(x))(\langle 1, z \rangle \in \mathsf{snd}(y)) \land (\forall w \in \mathsf{snd}(y))(\mathsf{snd}(w) \subseteq \mathsf{snd}(x)).$$
These definitions yield an interpretation, $I$, of an $\mathcal{L}_{\mathsf{P}}^{*}$-structure that is summarised in Table \ref{tab:table1} that extends the table in \cite[p. 373]{bar75}:
\begin{table}[H]
  \begin{center}
    \caption{The interpretation $I$ of an $\mathcal{L}_{\mathsf{P}}^{*}$-structure in $\mathcal{M}$}
    \label{tab:table1}
    \begin{tabular}{c|c} 
      {\bf $\mathcal{L}^*_\mathsf{P}$ Symbol} & {\bf $\mathcal{L}$ expression under $I$}\\
      \hline\\
      $\forall x$ & $\forall x(\mathsf{N}(x) \lor \mathsf{Set}(x) \Rightarrow \cdots)$ \\
      $=$ & $=$\\
      $\mathsf{U}(x)$ & $\mathsf{N}(x)$\\
      $x \mathsf{E} y$ & $x \mathsf{E}^\prime y$\\
      $x \in y$ & $x \mathcal{E} y$\\
      $\mathsf{F}(x)$ & $\bar{\mathsf{F}}(x)$\\
	  $\mathsf{P}(x)$ & $\bar{\mathsf{P}}(x)$
    \end{tabular}
  \end{center}
\end{table}
In other words, $\mathfrak{A}_\mathcal{N}= \langle \mathcal{N};\mathsf{Set}^\mathcal{M}, \mathcal{E}^\mathcal{M}, \bar{\mathsf{F}}^\mathcal{M}, \bar{\mathsf{P}}^\mathcal{M} \rangle$, where $\mathcal{N}= \langle \mathsf{N}^\mathcal{M}, (\mathsf{E}^\prime)^\mathcal{M} \rangle$, is an $\mathcal{L}_{\mathsf{P}}^{*}$-structure. If $\phi$ is an $\mathcal{L}_{\mathsf{P}}^{*}$-formula, then we write $\phi^I$ for the translation of $\phi$ into an $\mathcal{L}$-formula of $\mathcal{M}$ described in Table 1. Note that the map $x \mapsto \langle 0, x \rangle$ is an isomorphism between $\mathcal{M}$ and $\mathcal{N}$. The following is the refinement of \cite[Appendix Lemma 3.2]{bar75} noted by \cite[Chapter 2]{res87}:

\begin{Theorems1} \label{Th:KPUInInterpretation}
$\mathfrak{A}_{\mathcal{N}} \models \mathsf{KPU}_{\mathbb{C}\mathsf{ov}}$. \Square
\end{Theorems1}

We now turn to showing that axioms and axiom schemes transfer from $\mathcal{M}$ to $\mathfrak{A}_{\mathcal{N}}$.

\begin{Lemma1} \label{Th:PowersetInInterpretation}
$\mathfrak{A}_\mathcal{N} \models \mathsf{Powerset}$.
\end{Lemma1}

\begin{proof}
Let $a$ be a set of $\mathfrak{A}_{\mathcal{N}}$. To see that $\bar{\mathsf{P}}(a)$ exists, note that $a= \langle 1, a_0 \rangle$ and $\bar{\mathsf{P}}(a)= \langle 1, X \rangle$ where $X= \{1\} \times \mathcal{P}(a_0)$. Therefore, the powerset axiom in $\mathcal{M}$ ensures that $\bar{\mathsf{P}}$ is total in $\mathfrak{A}_{\mathcal{N}}$. Now, let $b$ be a set of $\mathfrak{A}_\mathcal{N}$. Work inside $\mathcal{M}$. Now, $b= \langle 1, b_0 \rangle$. And,
$$b \mathcal{E} \bar{\mathsf{P}}(a) \textrm{ iff } b_0 \subseteq a_0,$$
$$\textrm{iff for all } x, \textrm{ if } x \mathcal{E} b, \textrm{ then } x\mathcal{E} a,$$
$$\textrm{iff } (b \subseteq a)^I.$$
Therefore, $\mathfrak{A}_{\mathcal{N}}$ satisfies $\mathsf{Powerset}$.
\Square
\end{proof}

\begin{Lemma1} \label{Th:BoundedInPowersetOperationIsDelta1P}
Let $\phi(\vec{x})$ be a $\Delta_0(\mathcal{L}_{\mathsf{P}}^{*})$-formula. Then $\phi^I(\vec{x})$ is equivalent to a $\Delta_1^\mathcal{P}$-formula in $\mathcal{M}$.
\end{Lemma1}

\begin{proof}
We prove this lemma by induction on the complexity of $\phi^I$. Note that, by the above observations, $\mathsf{N}(x)$ and $x \mathsf{E}^\prime y$ can be written as $\Delta_0$-formulae. Moreover, $y= \bar{\mathsf{F}}(x)$ is equivalent to a $\Delta_0$-formula, and $y= \bar{\mathsf{P}}(x)$ is equivalent to a $\Delta_0^\mathcal{P}$-formula. Now, $y \mathcal{E} \bar{\mathsf{F}}(x)$ iff
$$\mathsf{fst}(y)= 0 \land \mathsf{snd}(y) \in \mathsf{snd}(x),$$
which is $\Delta_0$. Similarly, $y \mathcal{E} \bar{\mathsf{P}}(x)$ iff
$$\mathsf{fst}(y)= 1 \land \mathsf{snd}(y) \subseteq \mathsf{snd}(x),$$
which is also $\Delta_0$. Now, suppose that $t(x)$ is an $\mathcal{L}_{\mathsf{P}}^{*}$-term and both $y=t^I(x)$ and $y \mathcal{E} t^I(x)$ are $\Delta_1^\mathcal{P}$ in $\mathcal{M}$. Now, $y= \bar{\mathsf{P}}(t^I(x))$
$$\textrm{iff~~}\exists w (w = t^I(x) \land y= \bar{\mathsf{P}}(w)),$$
$$\textrm{iff~~} \forall w(w= t^I(x) \Rightarrow y= \bar{\mathsf{P}}(w)).$$
Similarly, $y \mathcal{E} \bar{\mathsf{P}}(t^I(x))$
$$\textrm {iff~~}\exists w ( w = t^I(x) \land y \mathcal{E} \bar{\mathsf{P}}(w)),$$
$$\textrm{iff~~} \forall w( w= t^I(x) \Rightarrow y \mathcal{E} \bar{\mathsf{P}}(w)).$$
Therefore, both $y= \bar{\mathsf{P}}(t^I(x))$ and $y \mathcal{E} \bar{\mathsf{P}}(t^I(x))$ are $\Delta_1^\mathcal{P}$ in $\mathcal{M}$. Now, $y= \bar{\mathsf{F}}(t^I(x))$
$$\textrm{iff~~} \exists w(w= t^I(x) \land y=\bar{\mathsf{F}}(w)),$$
$$\textrm{iff~~} \forall w(w= t^I(x) \Rightarrow y= \bar{\mathsf{F}}(w)).$$
And, $y \mathcal{E} \bar{\mathsf{F}}(t^I(x))$
$$\textrm{iff~~}\exists w ( w = t^I(x) \land y \mathcal{E} \bar{\mathsf{F}}(w)),$$
$$\textrm{iff~~} \forall w( w= t^I(x) \Rightarrow y \mathcal{E} \bar{\mathsf{F}}(w)).$$
Since $\bar{\mathsf{F}}$ and $\bar{\mathsf{P}}$ are both unary functions, this shows that for every $\mathcal{L}_{\mathsf{P}}^{*}$-term $t(x)$, both $y= t^I(x)$ and $y \mathcal{E} t^I(x)$ are $\Delta_1^\mathcal{P}$ in $\mathcal{M}$. Finally, we need an induction step that allows us to deal with bounded quantification. Let $\psi(x_0, \ldots, x_{n-1})$ be an $\mathcal{L}_{\mathsf{P}}^{*}$-formula such that $\psi^I(x_0, \ldots, x_{n-1})$ is $\Delta_1^\mathcal{P}$ in $\mathcal{M}$. Now, $(\exists x_0 \mathcal{E} x_n) \psi^I(x_0, \ldots, x_{n-1})$
$$\textrm{iff } (\exists x_0 \in \mathsf{snd}(x_n))\psi^I(x_0, \ldots, x_{n-1}).$$
Therefore, $(\exists x_0 \mathcal{E} x_n) \psi^I(x_0, \ldots, x_{n-1})= ((\exists x_0 \in x_n) \psi(x_0, \ldots, x_{n-1}))^I$ is $\Delta_1^{\mathcal{P}}$ in $\mathcal{M}$. Let $t(x)$ be an $\mathcal{L}_{\mathsf{P}}^{*}$-term. Now, $(\exists x_0 \mathcal{E} t^I(x_n)) \psi^I(x_0, \ldots, x_{n-1})$
$$\textrm{iff } \exists w(w = t^I(x_n) \land (\exists x_0 \in \mathsf{snd}(w)) \psi^I(x_0, \ldots, x_{n-1})),$$
$$\textrm{iff } \forall w (w= t^I(x_n) \Rightarrow (\exists x_0 \in \mathsf{snd}(w)) \psi^I(x_0, \ldots, x_{n-1})).$$
Therefore $(\exists x_0 \mathcal{E} t^I(x_n)) \psi^I(x_0, \ldots, x_{n-1})= ((\exists x_0 \in t(x_n)) \psi(x_0, \ldots, x_{n-1}))^I$ is $\Delta_1^\mathcal{P}$ in $\mathcal{M}$. The Lemma now follows by induction.
\Square
\end{proof}

\begin{Lemma1} \label{Th:SeparationInInterpretation}
$\mathfrak{A}_{\mathcal{N}} \models \Delta_0(\mathcal{L}_{\mathsf{P}}^{*})\text{-}\mathsf{Separation}$.
\end{Lemma1}

\begin{proof}
Let $\phi(x, \vec{z})$ be a $\Delta_0(\mathcal{L}_{\mathsf{P}}^{*})$-formula, $\vec{v}$ be sets and/or urelements of $\mathfrak{A}_{\mathcal{N}}$ and $a$ a set of $\mathfrak{A}_{\mathcal{N}}$. Work inside $\mathcal{M}$. Now, $a= \langle 1, a_0 \rangle$. Let
$$b_0= \{x \in a_0 \mid \phi^I(x, \vec{v})\},$$
which is a set by $\Delta_1^\mathcal{P}$-\textsf{Separation}. Let $b= \langle 1, b_0 \rangle$. Therefore, for all $x$ such that $\mathsf{Set}(x)$,
$$x \mathcal{E} b \textrm{ iff } x \mathcal{E} a \land \phi^I(x, \vec{v}).$$
Therefore, $\mathfrak{A}_{\mathcal{N}}$ satisfies $\Delta_0(\mathcal{L}_{\mathsf{P}}^{*})$-\textsf{Separation}.
\Square
\end{proof}

\begin{Lemma1} \label{Th:CollectionInInterpretation}
$\mathfrak{A}_\mathcal{N} \models \Delta_0(\mathcal{L}_{\mathsf{P}}^{*})\text{-}\mathsf{Collection}$.
\end{Lemma1}

\begin{proof}
Let $\phi(x, y, \vec{z})$ be a $\Delta_0(\mathcal{L}_{\mathsf{P}}^{*})$-formula. Let $\vec{v}$ be a sequence of sets and/or urelements of $\mathfrak{A}_\mathcal{N}$ and let $a$ be a set of $\mathfrak{A}_\mathcal{N}$ such that
$$\mathfrak{A}_\mathcal{N} \models (\forall x \in a) \exists y \phi(x, y, \vec{v}).$$
Work inside $\mathcal{M}$. Since $a$ is a set of $\mathfrak{A}_{\mathcal{N}}$, $a= \langle 1, a_0 \rangle$. We have
$$(\forall x \mathcal{E} a) \exists y ((\mathsf{N}(y) \lor \mathsf{Set}(y)) \land \phi^I(x, y, \vec{v})).$$
And,
$$(\forall x \in a_0) \exists y ((\mathsf{N}(y) \lor \mathsf{Set}(y)) \land \phi^I(x, y, \vec{v})).$$
So, since $(\mathsf{N}(y) \lor \mathsf{Set}(y)) \land \phi^I(x, y, \vec{v})$ is equivalent to a $\Sigma_1^\mathcal{P}$-formula, we can apply $\Delta_0^\mathcal{P}$-\textsf{Collection} to obtain $b$ such that
$$(\forall x \in a_0) (\exists y \in b) ((\mathsf{N}(y) \lor \mathsf{Set}(y)) \land \phi^I(x, y, \vec{v}))^{(b)}.$$
Let $b_0= \{y \in b \mid (\mathsf{N}(y) \lor \mathsf{Set}(y))^{(b)}\}$, which is a set by $\Delta_0$-\textsf{Separation}. Let $b_1= \langle 1, b_0 \rangle$. Therefore $\mathsf{Set}(b_1)$ and
$$(\forall x \mathcal{E} a)(\exists y \mathcal{E} b_1) \phi^I(x, y, \vec{v}).$$
So,
$$\mathfrak{A}_\mathcal{N} \models (\forall x \in a) (\exists y \in b_1) \phi(x, y, \vec{v}).$$
This shows that $\mathfrak{A}_{\mathcal{N}}$ satisfies $\Delta_0(\mathcal{L}_{\mathsf{P}}^{*})$-\textsf{Collection}.
\Square
\end{proof}

\begin{Lemma1} \label{Th:FoundationInterpretation}
$\mathfrak{A}_\mathcal{N} \models \Sigma_1(\mathcal{L}_{\mathsf{P}}^{*})\text{-}\mathsf{Foundation}$.
\end{Lemma1}

\begin{proof}
Let $\phi(x, \vec{z})$ be a $\Sigma_1(\mathcal{L}_{\mathsf{P}}^{*})$-formula. Let $\vec{v}$ be a sequence of sets and/or urelements be such that
$$\{ x \in \mathfrak{A}_{\mathcal{N}} \mid \mathfrak{A}_{\mathcal{N}} \models \phi(x, \vec{v}) \} \textrm{ is nonempty}.$$
Work inside $\mathcal{M}$. Consider $\theta(\alpha, \vec{z})$ defined by
$$(\alpha \textrm{ is an ordinal}) \land \exists x((\mathsf{Set}(x) \lor \mathsf{N}(x))\land \rho(x)=\alpha \land \phi^I(x, \vec{z})).$$
Note that $\theta(\alpha, \vec{z})$ is equivalent to a $\Sigma_1^\mathcal{P}$-formula. Therefore, using $\Sigma_1^\mathcal{P}$-\textsf{Foundation}, let $\beta$ be an $\in$-least element of
$$\{\alpha \in M \mid \mathcal{M} \models \theta(\alpha, \vec{v})\}.$$
Let $y$ be such that $(\mathsf{N}(y) \lor \mathsf{Set}(y))$, $\rho(y)= \beta$ and $\phi^I(y, \vec{v})$. Note that if $x \mathcal{E} y$, then $\rho(x) < \rho(y)$. Therefore $y$ is an $\mathcal{E}$-least element of
$$\{ x \in \mathfrak{A}_{\mathcal{N}} \mid \mathfrak{A}_{\mathcal{N}} \models \phi(x, \vec{v}) \}.$$
\Square
\end{proof}

The following combines \cite[II.8.4]{bar75} with the characterisation of $\mathbb{C}\mathsf{ov}_\mathcal{M}$ proved in \cite[Appendix Section 3]{bar75}:

\begin{Theorems1} \label{Th:BarwiseCharacterisationOfCover}
The $\mathcal{L}^*$-reduct of $\mathrm{WF}(\mathfrak{A}_\mathcal{N})$, $\mathrm{WF}^-(\mathfrak{A}_\mathcal{N})= \langle \mathcal{N}; \mathsf{WF}(\mathsf{Set}^\mathcal{M}), \mathcal{E}^\mathcal{M}, \bar{\mathsf{F}}^\mathcal{M} \rangle$ is an admissible set covering $\mathcal{N}$ that is isomorphic to $\mathbb{C}\mathsf{ov}_\mathcal{M}$. \Square
\end{Theorems1}

We now turn to extending this result to show that $\mathrm{WF}(\mathfrak{A}_\mathcal{N})$ is a power admissible set covering $\mathcal{N}$ and therefore the least power admissible set covering $\mathcal{N}$.

\begin{Theorems1} \label{Th:WellFoundedPartIsPowerAdmissibleCover}
The structure $\mathrm{WF}(\mathfrak{A}_\mathcal{N})= \langle \mathcal{N}; \mathrm{WF}(\mathsf{Set}^\mathcal{M}), \mathcal{E}^\mathcal{M}, \bar{\mathsf{F}}^\mathcal{M}, \bar{\mathsf{P}}^\mathcal{M} \rangle$ is a power admissible set covering $\mathcal{N}$. Moreover, $\mathrm{WF}(\mathfrak{A}_\mathcal{N})$ is isomorphic to $\mathbb{C}\mathsf{ov}^\mathsf{P}_\mathcal{M}$.
\end{Theorems1}

\begin{proof}
Note that it follows immediately from Theorem \ref{Th:BarwiseCharacterisationOfCover} that\\ $\mathrm{WF}(\mathfrak{A}_\mathcal{N})= \langle \mathcal{N}; \mathrm{WF}(\mathsf{Set}^\mathcal{M}), \mathcal{E}^\mathcal{M}, \bar{\mathsf{F}}^\mathcal{M}, \bar{\mathsf{P}}^\mathcal{M} \rangle$ satisfies all of the axioms of $\mathsf{KPU}_{\mathbb{C}\mathsf{ov}}$ plus full \textsf{Foundation}. The fact that $\mathrm{WF}(\mathfrak{A}_{\mathcal{N}}) \subseteq_e^\mathcal{P} \mathfrak{A}_{\mathcal{N}}$ implies that $\mathrm{WF}(\mathfrak{A}_{\mathcal{N}})$ satisfies \textsf{Powerset} and $\Delta_0(\mathcal{L}_{\mathsf{P}}^{*})$-\textsf{Separation}. To show that $\mathrm{WF}(\mathfrak{A}_{\mathcal{N}})$ satisfies $\Delta_0(\mathcal{L}_{\mathsf{P}}^{*})$-\textsf{Collection}, let $\phi(x, y, \vec{z})$ be a $\Delta_0(\mathcal{L}_{\mathsf{P}}^{*})$-formula. Let $\vec{v}$ be sets and/or urelements of $\mathsf{WF}(\mathfrak{A}_{\mathcal{N}})$ and let $a$ be a set of $\mathrm{WF}(\mathfrak{A}_{\mathcal{N}})$ such that
$$\mathrm{WF}(\mathfrak{A}_{\mathcal{N}}) \models (\forall x \in a) \exists y \phi(x, y, \vec{v}).$$
Consider the formula $\theta(\beta, \vec{z})$ defined by
$$(\beta \textrm{ is an ordinal})\land (\forall x \in a)(\exists \alpha \in \beta) \exists y (\rho(y)=\alpha \land \phi(x, y, \vec{z}).$$
Since $\mathrm{WF}(\mathfrak{A}_{\mathcal{N}}) \subseteq_e^\mathcal{P} \mathfrak{A}_{\mathcal{N}}$, if $\beta$ is a nonstandard ordinal of $\mathfrak{A}_{\mathcal{N}}$, then $\mathfrak{A}_\mathcal{N} \models \theta(\beta, \vec{v})$. Using $\Delta_0(\mathcal{L}_{\mathsf{P}}^{*})$-\textsf{Collection}, $\theta(\beta, \vec{z})$ is equivalent to a $\Sigma_1(\mathcal{L}_{\mathsf{P}}^{*})$-formula in $\mathfrak{A}_{\mathcal{N}}$. Therefore, by $\Sigma_1(\mathcal{L}_{\mathsf{P}}^{*})$-\textsf{Foundation}, $\{\beta \mid \mathfrak{A}_{\mathcal{N}} \models \theta(\beta, \vec{v})\}$ has a least element $\gamma$. Note that $\gamma$ is an ordinal of $\mathrm{WF}(\mathfrak{A}_\mathcal{N})$. Now, consider the formula $\psi(x, y, \vec{z}, \gamma)$ defined by
$$\phi(x, y, \vec{z}) \land (\rho(y) < \gamma).$$
Note that
$$\mathfrak{A}_\mathcal{N} \models (\forall x \in a) \exists y \psi(x, y, \vec{v}, \gamma).$$
Using $\Delta_0(\mathcal{L}_{\mathsf{P}}^{*})$-\textsf{Collection} in $\mathfrak{A}_\mathcal{N}$, there exists a set $b$ of $\mathfrak{A}_{\mathcal{N}}$ such that
$$\mathfrak{A}_\mathcal{N} \models (\forall x \in a)(\exists y \in b) \psi(x, y, \vec{v}, \gamma).$$
Let $c= \{x \in b \mid \rho(x) < \gamma\}$, which is a set in $\mathfrak{A}_\mathcal{N}$ by $\Delta_1(\mathcal{L}_{\mathsf{P}}^{*})$-\textsf{Separation}. Now, $c$ is a set of $\mathrm{WF}(\mathfrak{A}_\mathcal{N})$ and
$$\mathrm{WF}(\mathfrak{A}_\mathcal{N}) \models (\forall x \in a)(\exists y \in c) \phi(x, y, \vec{v}).$$
Therefore, $\mathsf{WF}(\mathfrak{A}_\mathcal{N})$ satisfies $\Delta_0(\mathcal{L}_{\mathsf{P}}^{*})$-\textsf{Collection}. And so, $\mathrm{WF}(\mathfrak{A}_{\mathcal{N}})$ is a power admissible set covering $\mathcal{N}$. Finally, since the $\mathcal{L}^*$-reduct of $\mathrm{WF}(\mathfrak{A}_{\mathcal{N}})$ is isomorphic to $\mathbb{C}\mathsf{ov}_\mathcal{M}$, $\mathrm{WF}(\mathfrak{A}_{\mathcal{N}})$ is isomorphic to $\mathbb{C}\mathsf{ov}^\mathsf{P}_\mathcal{M}$.
\Square
\end{proof}

The following theorem summarises the analysis undertaken in this section:

\begin{Theorems1} \label{Th:KeyResultOnCover}
If $\mathcal{M} \models \mathsf{KP}+\mathsf{Powerset}+\Delta_0^\mathcal{P}\text{-}\mathsf{Collection}+\Sigma_1^\mathcal{P}\text{-}\mathsf{Foundation}$, then there is an interpretation of $\mathsf{P}$ in $\mathbb{C}\mathsf{ov}_\mathcal{M}$ that yields the power admissible set $\mathbb{C}\mathsf{ov}^\mathsf{P}_\mathcal{M}$. \Square
\end{Theorems1}

This yields a version of \cite[Corollary 2.4.]{bar75} that will be useful for the compactness arguments in the next section.


\begin{Theorems1} \label{Th:FirstKeyLemmaForCompactness}
Let $\mathcal{M}= \langle M, \mathsf{E}^\mathcal{M} \rangle \models \mathsf{KP}+\mathsf{Powerset}+\Delta_0^\mathcal{P}\text{-}\mathsf{Collection}+\Sigma_1^\mathcal{P}\text{-}\mathsf{Foundation}$. For all $A \subseteq M$, there exists $a \in M$ such that $a^*=A$ if and only if $A \in \mathbb{C}\mathsf{ov}^\mathsf{P}_\mathcal{M}$. \Square
\end{Theorems1}

\section[Powerset-preserving and rank extensions]{End extension results}

In this section we use the Barwise Compactness Theorem for $\mathcal{L}^{\mathsf{ee}}_{\mathbb{C}\mathsf{ov}^\mathsf{P}_\mathcal{M}}$ to show that every countable model of $\mathsf{KP}+\mathsf{powerset}+\Delta_0^\mathcal{P}\text{-}\mathsf{Collection}+\Sigma_1^\mathcal{P}\text{-}\mathsf{Foundation}$ has a powerset-preserving end extension.

The following is an immediate consequence of Theorem \ref{Th:FirstKeyLemmaForCompactness}:

\begin{Lemma1} \label{Th:SecondKeyLemmaForCompactness}
Let $\mathcal{M}= \langle M, \mathsf{E}^\mathcal{M} \rangle\models \mathsf{KP}+\mathsf{Powerset}+\Delta_0^\mathcal{P}\text{-}\mathsf{Collection}+\Sigma_1^\mathcal{P}\text{-}\mathsf{Foundation}$, and let $T_0$ be an $\mathcal{L}^{\mathsf{ee}}_{\mathbb{C}\mathsf{ov}^\mathsf{P}_\mathcal{M}}$-theory. If $T_0 \in \mathbb{C}\mathsf{ov}^\mathsf{P}_\mathcal{M}$, then there exists $b\in M$ such that
$$b^*= \{ a \in M \mid \bar{a} \textrm{ is mentioned in } T_0\}.$$
\Square
\end{Lemma1}

The next result expands on comments made in \cite[p. 637]{bar75} and connects definability in $\mathcal{M}$ to definability in $\mathbb{C}\mathsf{ov}^\mathsf{P}_\mathcal{M}$.

\begin{Lemma1} \label{Th:BoundInMTranslatesToBoundedInCov}
Let $\mathcal{M}= \langle M, \mathsf{E}^\mathcal{M} \rangle \models \mathsf{KP}+\mathsf{Powerset}+\Delta_0^\mathcal{P}\text{-}\mathsf{Collection}+\Sigma_1^{\mathcal{P}}\text{-}\mathsf{Foundation}$, and let $\phi(\vec{z})$ be a $\Delta_0^\mathcal{P}$-formula. Then there exists a formula $\hat{\phi}(\vec{z})$ that is $\Delta_1(\mathcal{L}^*_\mathsf{P})$ in the theory $\mathsf{KPU}^{\mathcal{P}}_{\mathbb{C}\mathsf{ov}}$ such that for all $\vec{z} \in M$,
$$\mathcal{M} \models \phi(\vec{z}) \text{~~iff~~} \mathbb{C}\mathsf{ov}^\mathsf{P}_\mathcal{M} \models \hat{\phi}(\vec{z}).$$
\end{Lemma1}

\begin{proof}
Let $\phi(\vec{z})$ be a $\Delta_0$-formula. We prove the lemma by structural induction on the complexity of $\phi$. Without loss of generality we can assume that the only connectives of propositional logic appearing in $\phi$ are $\neg$ and $\lor$. If $\phi(z_1, z_2)$ is $z_1 \in z_2$, then let $\hat{\phi}(z_1, z_2)$ be the $\Delta_0(\mathcal{L}^*_\mathsf{P})$-formula $z_1 \mathsf{E} z_2$. Therefore, for all $z_1, z_2 \in M$,
$$\mathcal{M} \models \phi(z_1, z_2) \textrm{~~iff~~} \mathbb{C}\mathsf{ov}^\mathsf{P}_\mathcal{M} \models \hat{\phi}(z_1, z_2).$$
If $\phi(\vec{z})$ is $\neg \psi(\vec{z})$ and the lemma holds for $\psi(\vec{z})$, then let $\hat{\phi}(\vec{z})= \neg \hat{\psi}(\vec{z})$. So, $\hat{\phi}(\vec{z})$ is $\Delta_1(\mathcal{L}^*_\mathsf{P})$ in the theory $\mathsf{KPU}^{\mathcal{P}}_{\mathbb{C}\mathsf{ov}}$ and for all $\vec{z} \in M$,
$$\mathcal{M} \models \phi(\vec{z}) \textrm{~~iff~~} \mathbb{C}\mathsf{ov}^\mathsf{P}_\mathcal{M} \models \hat{\phi}(\vec{z}).$$
Suppose that $\phi(\vec{z})$ is $\psi_1(\vec{z}) \lor \psi_2(\vec{z})$ and the lemma holds for $\psi_1(\vec{z})$ and $\psi_2(\vec{z})$. Let $\hat{\phi}(\vec{z})$ be $\hat{\psi_1}(\vec{z}) \lor \hat{\psi_2}(\vec{z})$. Therefore, $\hat{\phi}(\vec{z})$ is $\Delta_1(\mathcal{L}^*_\mathsf{P})$ in the theory $\mathsf{KPU}^{\mathcal{P}}_{\mathbb{C}\mathsf{ov}}$ and for all $\vec{z} \in M$,
$$\mathcal{M} \models \phi(\vec{z}) \textrm{~~iff~~ } \mathbb{C}\mathsf{ov}^\mathsf{P}_\mathcal{M} \models \hat{\phi}(\vec{z}).$$
Suppose $\phi(y, \vec{z})$ is $(\mathcal{Q} x \in y)\psi(x, y, \vec{z})$, where $\mathcal{Q} \in \{\exists, \forall\}$, and the lemma holds for $\psi(x, y, \vec{z})$. Let $\hat{\phi}(y, \vec{z})$ be $(\mathcal{Q} x \in \mathsf{F}(y))\hat{\psi}(x, y, \vec{z})$. So, $\hat{\phi}(\vec{z})$ is $\Delta_1(\mathcal{L}^*_\mathsf{P})$ in the theory $\mathsf{KPU}^{\mathcal{P}}_{\mathbb{C}\mathsf{ov}}$. Since $\mathbb{C}\mathsf{ov}^\mathsf{P}_\mathcal{M}$ satisfies $(\dagger)$, for all $y, \vec{z} \in M$,
$$\mathcal{M} \models \phi(y, \vec{z}) \textrm{~~iff~~} \mathcal{M} \models (\mathcal{Q} x \in y) \psi(x, y, \vec{z})$$
$$\textrm{~~iff~~} \mathbb{C}\mathsf{ov}^\mathsf{P}_\mathcal{M} \models (\mathcal{Q} x \in \mathsf{F}(y))\hat{\psi}(x, y, \vec{z})$$
$$\textrm{~~iff~~}\mathbb{C}\mathsf{ov}^\mathsf{P}_\mathcal{M} \models \hat{\phi}(y, \vec{z}).$$
Suppose that $\phi(y, \vec{z})$ is $(\mathcal{Q} x \subseteq y) \psi(x, y, \vec{z})$, where $\mathcal{Q} \in \{\exists, \forall\}$, and the lemma holds for $\psi(x, y, \vec{z})$. Let $\hat{\phi}(y, \vec{z})$ be $(\mathcal{Q} x \in \mathsf{P}(\mathsf{F}(y)))\exists p(\mathsf{F}(p)= x \land \hat{\psi}(p, y, \vec{z}))$. Note that, in the theory $\mathsf{KPU}^{\mathcal{P}}_{\mathbb{C}\mathsf{ov}}$, for all urelements $\vec{z}$,
$$(\mathcal{Q} x \in \mathsf{P}(\mathsf{F}(y)))~\exists p(\mathsf{F}(p)= x ~\land ~ \hat{\psi}(p, y, \vec{z})) \iff (\mathcal{Q} x \in \mathsf{P}(\mathsf{F}(y)))~\forall p(\mathsf{F}(p)= x \Rightarrow \hat{\psi}(p, y, \vec{z})).$$
Therefore, in the theory $\mathsf{KPU}^{\mathcal{P}}_{\mathbb{C}\mathsf{ov}}$, $\hat{\phi}(y, \vec{z})$ is $\Delta_1(\mathcal{L}^*_\mathsf{P})$. Moreover, by Theorem \ref{Th:FirstKeyLemmaForCompactness} and $(\dagger)$, for all $\vec{z} \in M$,
$$\mathcal{M} \models \hat{\phi}(y, \vec{z}) \textrm{~~iff~~} \mathcal{M} \models (\mathcal{Q} x \subseteq y) \psi(x, y, \vec{z})$$
$$\textrm{~~iff~~} (\mathcal{Q} x \in \mathsf{P}(\mathsf{F}(y)))~\exists p(\mathsf{F}(p)= x \land \hat{\psi}(p, y, \vec{z}))$$
$$\textrm{~~iff~~}\mathbb{C}\mathsf{ov}^\mathsf{P}_\mathcal{M} \models \hat{\phi}(y, \vec{z}).$$
Therefore, the lemma follows by induction.
\Square
\end{proof}

We are now able to use the machinery we have developed to establish the following result.



\begin{Theorems1} \label{Th:MainEndExtensionResult}
Let $S$ be a recursively enumerable $\mathcal{L}$-theory such that $$S \vdash \mathsf{KP}+\mathsf{Powerset}+\Delta_0^\mathcal{P}\text{-}\mathsf{Collection}+\Sigma_1^\mathcal{P}\text{-}\mathsf{Foundation},$$and let $\mathcal{M}$ be a countable model of $S$. Then there exists an $\mathcal{L}$-structure $\mathcal{N}$ such that $\mathcal{M} \subseteq_e^\mathcal{P} \mathcal{N}\models S$, and for some $d \in N$, and for all $x \in M$, $\mathcal{N} \models (x \in d)$.
\end{Theorems1}



\begin{proof}
Let $T$ be the $\mathcal{L}^{\mathsf{ee}}_{\mathbb{C}\mathsf{ov}^\mathsf{P}_\mathcal{M}}$-theory that contains:
\begin{itemize}
\item $S$;
\item for all $a, b \in M$ with $\mathcal{M} \models (a \in b)$, $\bar{a} \in \bar{b}$;
\item for all $a \in M$,
$$\forall x\left( x \in \bar{a} \iff \bigvee_{b \in a} (x= \bar{b}) \right);$$
\item for all $a \in M$,
$$\forall x\left( x \subseteq \bar{a} \iff \bigvee_{b \subseteq a} (x= \bar{b}) \right);$$
\item for all $a \in M$, $\bar{a} \in \mathbf{c}$.
\end{itemize}
Lemma \ref{Th:BoundInMTranslatesToBoundedInCov} shows that $T \subseteq \mathbb{C}\mathsf{ov}^\mathsf{P}_\mathcal{M}$ is $\Sigma_1(\mathcal{L}_{\mathsf{P}}^{*})$-definable over $\mathbb{C}\mathsf{ov}^\mathsf{P}_\mathcal{M}$. Let $T_0 \subseteq T$ be such that $T_0 \in \mathbb{C}\mathsf{ov}^\mathsf{P}_\mathcal{M}$. Using Lemma \ref{Th:SecondKeyLemmaForCompactness}, there exists $c \in M$ such that
$$c^*= \{a \in M \mid \bar{a} \textrm{ is mentioned in } T_0\}.$$
Therefore, by interpreting each $\bar{a}$ that is mentioned in $T_0$ by $a \in M$ and interpreting $\mathbf{c}$ by $c$, we can expand $\mathcal{M}$ to a model $\mathcal{M}^\prime$ that satisfies $T_0$. Therefore, by the Barwise Compactness Theorem, there exists $\mathcal{N} \models T$. It is straightforward to see that the $\mathcal{L}$-reduct of $\mathcal{N}$ is the desired extension of $\mathcal{M}$.
\Square
\end{proof}

We first apply this result to show that countable models of $\mathsf{KP}^\mathcal{P}$ have topless rank extensions that satisfy $\mathsf{KP}^\mathcal{P}$. This generalises \cite[Theorem 2.3]{fri73}, which shows that every countable transitive model of $\mathsf{KP}^\mathcal{P}$ has a topless rank extension that satisfies $\mathsf{KP}^\mathcal{P}$.

\begin{Theorems1} \label{Th:RankExtensionsCtbleTransModelsOfKPP}
(Friedman) Let $S$ be a recursively enumerable $\mathcal{L}$-theory such that $S \vdash \mathsf{KP}^\mathcal{P}$. If $\mathcal{M}$ is a countable transitive model of $S$, then there exists $\mathcal{N} \models S$ such that $\mathcal{M} \subseteq_{\mathsf{topless}}^{\mathsf{rk}} \mathcal{N}$. \Square
\end{Theorems1}

It follows from \cite[Theorem 4.8]{gor20} that every countable nonstandard model of $\mathsf{KP}^\mathcal{P}+\Sigma_1^\mathcal{P}\textsf{-Separation}$, $\mathcal{N}$, is isomorphic to substructure $\mathcal{M}$ of $\mathcal{N}$ such that $\mathcal{M} \subseteq_{\mathsf{topless}}^\mathsf{rk} \mathcal{N}$. We will make use of this result in the following form:

\begin{Theorems1} \label{Th:GorbowSelfEmbeddingTheorem}
(Gorbow) Let $\mathcal{M}$ be a countable nonstandard model of $\mathsf{KP}^\mathcal{P}+\Sigma_1^\mathcal{P}\text{-}\mathsf{Separation}$. Then there exists $\mathcal{N} \equiv \mathcal{M}$ such that $\mathcal{M} \subseteq_{\mathsf{topless}}^{\mathsf{rk}} \mathcal{N}$. \Square
\end{Theorems1}

We next note that if a model of $\mathsf{KP}^\mathcal{P}$ has a blunt rank extension, then that model must satisfy the full scheme of separation.

\begin{Lemma1} \label{Th:BluntExtensionImpliesSeparation}
Let $\mathcal{M}= \langle M, \mathsf{E}^\mathcal{M} \rangle$ and $\mathcal{N}= \langle N, \mathsf{E}^\mathcal{N} \rangle$ be models of $\mathsf{KP}^\mathcal{P}$. If $\mathcal{M} \subseteq_{\mathsf{blunt}}^{\mathsf{rk}} \mathcal{N}$, then $\mathcal{M} \models \Pi_\infty\text{-}\mathsf{Separation}$.
\end{Lemma1}

\begin{proof}
Assume that $M \subseteq N$, $\mathsf{E}^\mathcal{M}= \mathsf{E}^\mathcal{N} \upharpoonright M$ and $\mathcal{M} \subseteq_{\mathsf{blunt}}^{\mathsf{rk}} \mathcal{N}$. Let $c \in N$ be such that $c^* \subseteq M$ and $c \notin M$. Working inside $\mathcal{N}$, let $\alpha= \rho(c)$. Therefore, since $\mathcal{M} \subseteq_{\mathsf{blunt}}^{\mathsf{rk}} \mathcal{N}$,
$$\begin{array}{lll}
x \in (V_\alpha^\mathcal{N})^* & \textrm{if and only if} & \mathcal{N} \models (\rho(x) < \alpha)\\
& \textrm{if and only if} & \mathcal{N} \models (\exists y \in c)(\rho(x) \leq \rho(y))\\
& \textrm{if and only if} & x \in M.
\end{array}$$
So, $M= (V_\alpha^\mathcal{N})^*$ and every instance of $\Pi_\infty$-\textsf{Separation} in $\mathcal{M}$ can be reduced to an instance of $\Delta_0$-\textsf{Separation} in $\mathcal{N}$ and, since $\mathcal{M} \subseteq_e^\mathcal{P} \mathcal{N}$, the resulting set will be in $\mathcal{M}$. Therefore, $\mathcal{M} \models \Pi_\infty\text{-}\mathsf{Separation}$.
\Square
\end{proof}



\begin{Theorems1} \label{Th:ToplessEndExtensionsOfKPP}
Let $S$ be a recursively enumerable $\mathcal{L}$-theory such $S \vdash \mathsf{KP}^\mathcal{P}$. If $\mathcal{M}$ is a countable model of $S$, then there exists $\mathcal{N} \models S$ such that $\mathcal{M} \subseteq_{\mathsf{topless}}^{\mathsf{rk}} \mathcal{N}$.
\end{Theorems1}

\begin{proof}
Let $\mathcal{M}$ be a countable model of $S$. If $\mathcal{M}$ is well-founded, then $\mathcal{M}$ is isomorphic to a transitive model of $\mathsf{KP}^\mathcal{P}$ and we can use Theorem \ref{Th:RankExtensionsCtbleTransModelsOfKPP} to find an $\mathcal{L}$-structure $\mathcal{N} \models S$ such that $\mathcal{M} \subseteq_{\mathsf{topless}}^{\mathsf{rk}} \mathcal{N}$. Therefore, assume that $\mathcal{M}$ is nonstandard. By Corollary \ref{Th:KPPProvesSigma1PFoundation}, $\mathcal{M}$ satisfies $\mathsf{KP}+\mathsf{Powerset}+\Delta_0^\mathcal{P}\text{-}\mathsf{Collection}+\Sigma_1^\mathcal{P}\text{-}\mathsf{Foundation}$. Therefore, using Theorem \ref{Th:MainEndExtensionResult}, we can find an $\mathcal{L}$-structure $\mathcal{N} \models S$ such that $M \neq N$ and $\mathcal{M} \subseteq_e^\mathcal{P} \mathcal{N}$. So, by Lemma \ref{Th:PowersetPreservingEndExtensionsOfKPPAreRankExtensions}, $\mathcal{M} \subseteq_e^{\mathsf{rk}} \mathcal{N}$. If $\mathcal{M} \subseteq_{\mathsf{topless}}^{\mathsf{rk}} \mathcal{N}$, then we are done. Alternatively, if $\mathcal{M} \subseteq_{\mathsf{blunt}}^{\mathsf{rk}} \mathcal{N}$, then, by Lemma \ref{Th:BluntExtensionImpliesSeparation}, $\mathcal{M} \models \Pi_\infty\text{-}\mathsf{Separation}$. Therefore, since $\mathcal{M}$ is nonstandard, we can apply Theorem \ref{Th:GorbowSelfEmbeddingTheorem} to obtain an $\mathcal{L}$-structure $\mathcal{N}^\prime \equiv \mathcal{M}$ such that $\mathcal{M} \subseteq_{\mathsf{topless}}^{\mathsf{rk}} \mathcal{N}^\prime$.
\Square
\end{proof}

We now turn to showing that every countable model of $\mathsf{MOST}+\Pi_1\text{-}\mathsf{Collection}+\Sigma_1^\mathcal{P}\text{-}\mathsf{Foundation}$ has a topless powerset-preserving end extension that satisfies $\mathsf{MOST}+\Pi_1\text{-}\mathsf{Collection}+\Sigma_1^\mathcal{P}\text{-}\mathsf{Foundation}$. We first prove an analogue of Lemma \ref{Th:BluntExtensionImpliesSeparation} for models of $\mathsf{MOST}$.

\begin{Lemma1} \label{Th:BluntExtensionImpliesSeparationMOST}
Let $\mathcal{M}= \langle M, \mathsf{E}^\mathcal{M} \rangle$ and $\mathcal{N}= \langle N, \mathsf{E}^\mathcal{N} \rangle$ be models of $\mathsf{MOST}$. If $\mathcal{M} \subseteq_{\mathsf{blunt}}^{\mathcal{P}} \mathcal{N}$, then $\mathcal{M} \models \Pi_\infty\text{-}\mathsf{Separation}$.
\end{Lemma1}

\begin{proof}
Assume that $M \subseteq N$, $\mathsf{E}^\mathcal{M}=\mathsf{E}^\mathcal{N} \upharpoonright M$ and $\mathcal{M} \subseteq_{\mathsf{blunt}}^{\mathsf{rk}} \mathcal{N}$. Let $c \in N$ be such that $c^* \subseteq M$ and $c \notin M$. Work inside $\mathcal{N}$. Let $\kappa= |\mathsf{TC}(c)|$ and note that $\kappa \notin M$. Consider
$$A= \{ \lambda \in \kappa \mid (\exists y \in c)(\lambda= |\mathsf{TC}(y)|)\},$$
which is a set by $\Sigma_1$-\textsf{Separation}. Let $\mu= \sup A$ and note that $\mu$ is an initial ordinal. Work in the metatheory again. If $\mu \in M$, then so is $(\mu^+)^\mathcal{N}= (\mu^+)^\mathcal{M} \in M$. So, $H_{\mu^+}^\mathcal{N} \in M$ and $\mathcal{N} \models (c \subseteq H_{\mu^+})$. And, $c \in M$, which is a contradiction. Therefore $\mu \notin M$. Now,
$$\begin{array}{lll}
x \in (H_{\mu}^\mathcal{N})^* & \textrm{if and only if} & \mathcal{N} \models (|\mathsf{TC}(x)| < \mu)\\
& \textrm{if and only if} & \mathcal{N} \models (\exists y \in c)(|\mathsf{TC}(x)| < |\mathsf{TC}(y)|)\\
& \textrm{if and only if} & x \in M.
\end{array}$$
So, $M= (H_\mu^\mathcal{N})^*$ and every instance of $\Pi_\infty$-\textsf{Separation} in $\mathcal{M}$ can be reduced to an instance of $\Delta_0$-\textsf{Separation} in $\mathcal{N}$ and, since $\mathcal{M} \subseteq_e^\mathcal{P} \mathcal{N}$, the resulting set will be in $\mathcal{M}$. Therefore, $\mathcal{M} \models \Pi_\infty\text{-}\mathsf{Separation}$.
\Square
\end{proof}





\begin{Theorems1} \label{Th:ToplessPowersetPreservingEndExtensionOfMOST}
Let $S$ be an $\mathcal{L}$-theory such that
$S \vdash \mathsf{MOST}+\Pi_1\text{-}\mathsf{Collection}+\Sigma_1^\mathcal{P}\text{-}\mathsf{Foundation}$. If $\mathcal{M}$ is a countable model of $S$, then there exists  a model $\mathcal{N}$ such that $\mathcal{M} \subseteq_{\mathsf{topless}}^{\mathcal{P}} \mathcal{N}\models S$.
\end{Theorems1}

\begin{proof}
This can be proved using an identical argument to the proof of Theorem \ref{Th:ToplessEndExtensionsOfKPP} after observing that every transitive model of $\mathsf{MOST}+\Pi_1\text{-}\mathsf{Collection}+\Sigma_1^\mathcal{P}\text{-}\mathsf{Foundation}$ is a model of $\mathsf{KP}^\mathcal{P}$ and $\mathsf{KP}^\mathcal{P}+\Sigma_1^\mathcal{P}\text{-}\mathsf{Separation}$ is a subtheory of $\mathsf{MOST}+\Pi_1\text{-}\mathsf{Collection}+\Pi_\infty\text{-}\mathsf{Separation}$.
\Square
\end{proof}

The work \cite{ekm18} studies the class $\mathcal{C}$ of structures $\mathcal{I}_{\mathsf{fix}(j)}$ where $j: \mathcal{M} \longrightarrow \mathcal{M}$ is a nontrivial automorphism, $\mathcal{M}$ is an $\mathcal{L}$-structure that satisfies $\mathsf{MOST}$, $j$ fixes every point in $(\omega^\mathcal{M})^*$ and $\mathcal{I}_{\mathsf{fix}(j)}$ is the substructure of $\mathcal{M}$ that consists of elements $x$ of $\mathcal{M}$ such that $j$ fixes every point in $(\mathsf{TC}^\mathcal{M}(\{x\}))^*$. The results of \cite[Section 3]{ekm18} show that every structure in $\mathcal{C}$ satisfies $\mathsf{MOST}+\Pi_1\text{-}\mathsf{Collection}$. Conversely, \cite[Section 4]{ekm18} shows that a sufficient condition for a countable structure $\mathcal{M}$ that satisfies $\mathsf{MOST}+\Pi_1\text{-}\mathsf{Collection}$ to be in $\mathcal{C}$ is that there exists $\mathcal{M} \subseteq_{\mathsf{topless}}^\mathcal{P} \mathcal{N}$ such that $\mathcal{N}$ satisfies $\mathsf{MOST}+\Pi_1\text{-}\mathsf{Collection}$. Theorem \ref{Th:ToplessPowersetPreservingEndExtensionOfMOST} allows us to extend \cite[Theorem B]{ekm18} by showing that $\mathcal{C}$ contains all countable models of $\mathsf{MOST}+\Pi_1\text{-}\mathsf{Collection}+\Sigma_1^\mathcal{P}\text{-}\mathsf{Foundation}$.

\begin{Theorems1}\label{thm: aut}
Let $\mathcal{M}=\langle M, \mathsf{E}^\mathcal{M} \rangle$ be a countable model of $\mathsf{MOST}+\Pi_1\text{-}\mathsf{Collection}+\Sigma_1^\mathcal{P}\text{-}\mathsf{Foundation}$. Then there exists a model $\mathcal{N}=\langle N, \mathsf{E}^\mathcal{N} \rangle$ that satisfies $\mathsf{MOST}$ and a nontrivial automorphism $j:\mathcal{N} \longrightarrow \mathcal{N}$ such that $\mathcal{M} \cong \mathcal{I}_{\mathsf{fix}(j)}$, where $\mathcal{I}_{\mathsf{fix}(j)}$ is the substructure of $\mathcal{N}$ with underlying set
$$I_{\mathsf{fix}(j)}= \{x \in N \mid (\forall y \in (\mathsf{TC}^\mathcal{N}(\{x\}))^*)(j(y)=y)\}.$$
\Square
\end{Theorems1}

\noindent Combined with Corollary \ref{Th:Sigma1PFoundationInOmegaModelsMOST} and \cite[Theorem 5.6]{ekm18} this shows that the class $\mathcal{C}$ contains every countable recursively saturated model of $\mathsf{MOST}+\Pi_1\textrm{-\textsf{Collection}}$ and every countable $\omega$-standard model of $\mathsf{MOST}+\Pi_1\textrm{-\textsf{Collection}}$, providing a partial positive answer to Question 5.1 of \cite{ekm18}. A positive answer to the following question would positively answer Question 5.1 of \cite{ekm18}:

\begin{Quest1}
Does every countable $\omega$-nonstandard model of $\mathsf{MOST}+\Pi_1\text{-}\mathsf{Collection}$ have a topless powerset-preserving end extension that satisfies $\mathsf{MOST}+\Pi_1\text{-}\mathsf{Collection}$?
\end{Quest1}

Note that \cite[Theorem 5.6]{ekm18} shows that this question has a positive answer when the countable model is recursively saturated.

\bibliographystyle{alpha}
\bibliography{.}

\end{document}